\documentclass[a4paper,12pt]{article}
\usepackage{amsmath,amsfonts}
\usepackage{amsmath}
\usepackage{graphicx}
\usepackage{amsfonts}
\usepackage{amssymb}
\DeclareMathOperator{\expect}{{\mathbb E}}

\newcommand{\eps}{\varepsilon}

\makeatletter
\newcommand{\eqcolon}{\mathrel{\mathord{=}\raise.2\p@\hbox{:}}}
\newcommand{\coloneq}{\mathrel{\raise.2\p@\hbox{:}\mathord{=}}}
\makeatother

\newtheorem{lemma}{Lemma}
\newtheorem{proposition}{Proposition}
\newtheorem{corollary}{Corollary}
\newtheorem{theorem}{Theorem}
\newtheorem{remark}{Remark}

\begin{document}

\title{Statistics of a vortex filament model}
\author{{Franco Flandoli \& Massimiliano Gubinelli}\\{\small \textit{Dipartimento di Matematica Applicata}} \\{\small \textit{Universit\`{a} di Pisa, via Bonanno 25 bis}} \\{\small \textit{56126 Pisa, Italy}}\\{\small \textit{email:}\texttt{\{f.flandoli,m.gubinelli\}@dma.unipi.it}}}
\maketitle
\begin{abstract}
A random field composed by Poisson distributed Brownian vortex filaments is
constructed. The filament have a random thickness, length and intensity,
governed by a measure $\gamma$. Under appropriate assumptions on $\gamma$ we
compute the scaling law of the structure function and get the multifractal
scaling as a particular case.
\end{abstract}


\section{Introduction}

Isotropic homogeneous turbulence is phenomenologically described by several
theories, which usually give us the scaling properties of moments of velocity
increments. If $u(x)$ denotes the velocity field of the fluid and
$S_{p}\left(  \varepsilon\right)  $, the so called \textit{structure
function}, denotes the $p$-moment of the velocity increment over a distance
$\varepsilon$ (often only its longitudinal projection is considered), then one
expect a behavior of the form
\begin{equation}
S_{p}(\varepsilon) = \langle|u(x+\varepsilon)-u(x)|^{p} \rangle\sim
\varepsilon^{\zeta_{p}}. \label{scal}%
\end{equation}
Here, in our notations, $\varepsilon$ is not the dissipation energy, but just
the spatial scale parameter (see remark \ref{advert}). Let us recall two major
theories: the Kolmogorov-Obukov scaling law (K41) (see \cite{Kol}) says that
\[
\zeta_{2}=\frac{2}{3}%
\]
probably the best result compared with experiments; however the heuristic
basis of the theory also implies $\zeta_{p}=\frac{p}{3}$ which is not in
accordance with experiments. Intermittency corrections seem to be important
for larger $p$'s. A general theory which takes them into account is the
multifractal scaling theory of Parisi and Frisch \cite{Parisi-Frisch}, that
gives us $\zeta_{p}$ in the form of a Fenchel-Legendre transform:
\[
\zeta_{p}=\inf_{h\in I}[hp+3-D(h)].
\]
This theory is a sort of container, which includes for instance the striking
particular case of She and Leveque \cite{she}. We do not pretend to go further
in the explanation of this topic and address the reader to the monograph
\cite{frisch}.

The foundations of these theories, in particular of the multifractal one, are
usually mathematically poor, based mainly on very good intuition and a
suitable ``mental image'' (see the beginning of Chapter 7 of \cite{frisch}).
Essentially, the scaling properties of $S_{p}(\varepsilon)$ are given
\textit{a priori}, after an intuitive description of the mental image. The
velocity field of the fluid is not mathematically described or constructed,
but some crucial aspects of it are described only in plain words, and then
$S_{p}(\varepsilon)$ is given (or heuristically ``deduced'').

We do not pretend to remedy here to this extremely difficult problem, which
ultimately should start from a Navier-Stokes type model and the analysis of
its invariant measures.

The contribution of this paper is only to construct rigorously a random
velocity field which has two interesting properties: i) its realizations have
a geometry inspired by the pictures obtained by numerical simulations of
turbulent fluids; ii) the asymptotic as $\varepsilon\rightarrow0$ of
$S_{p}(\varepsilon)$ can be explicitly computed and the multifractal model is
recovered with a suitable choice of the measures defining the random field.
Its relation with Navier-Stokes models and their invariant measures is obscure
as well (a part from some vague conjectures, see \cite{Fl2}), so it is just
one small step beyond pure phenomenology of turbulence.

Concerning (i), the geometry of the field is that of a collection of vortex
filaments, as observed for instance by \cite{Vincent-Meneg} and many others.
The main proposal to model vortex filaments by paths of stochastic processes
came from A. Chorin, who made several considerations about their statistical
mechanics, see \cite{Ch}. The processes considered in \cite{Ch} are
self-avoiding walks, hence discrete. Continuous processes like Brownian
motion, geometrically more natural, have been considered by \cite{Gallavotti},
\cite{Lions-Majda}, \cite{Fl1}, \cite{FG1}, \cite{Nualart} and others. We do
not report here numerical results, but we have observed in simple simulations
that the vortex filaments of the present paper, with the tubular smoothening
due to the parameter $\ell$ (see below), have a shape that reminds very
strongly the simulations of \cite{Bell-Markus}.

Concerning (ii), we use stochastic analysis, properties of stochastic
integrals and ideas related to the theory of the Brownian sausage and
occupation measure. We do not know whether it is possible to reach so strict
estimates on $S_{p}(\varepsilon)$ as those proved here in the case when
stochastic processes are replaced by smooth curves. The power of stochastic
calculus seem to be important.

Ensembles of vortex structures with more stiff or artificial geometry have
been considered recently by \cite{Pullin} and \cite{HKambe}. They do not
stress the relation with multifractal models and part of their results are
numerical, but nevertheless they indicate that scaling laws can be obtained by
models based on many vortex structures. Probably a closer investigation of
simpler geometrical models like that ones, in spite of the less appealing
geometry of the objects, will be important to understand more about this
approach to turbulence theory. Let us also say that these authors introduced
that models also for numerical purposes, so the simplicity of the structures
has other important motivations.

Finally, let us remark about a difference with respect to the idea presented
in \cite{Ch} and also in \cite{Boyer}. There one put the attention on a single
vortex and try to relate statistical properties of the path of the process,
like Flory exponents of 3D self-avoiding walk, with scaling law of the
velocity field. Such an attempt is more intrinsic, in that it hopes to
associate turbulent scalings with relevant exponents known for processes. On
the contrary, here (and in \cite{Pullin} and \cite{HKambe}) we consider a
fluid made of a multitude of vortex structures and extract statistics from the
collective behavior. In fact, it this first work on the subject, we consider
independent vortex structures only, having in mind the Gibbs couplings of
\cite{FG1} as a second future step. Due to the independence, at the end again
it is just the single filament that determines, through the statistics of its
parameters, the properties of $S_{p}(\varepsilon)$. However, the
interpretation of the results and the conditions on the parameters are more in
the spirit of the classical ideas of K41 and its variants, where one thinks to
the 3D space more or less filled in by eddies or other structures. It is less
natural to interpret multifractality, for instance, on a single filament
(although certain numerical simulation on the evolution of a single filament
suggest that\ multifractality could arise on a single filament by a
non-uniform procedure of stretching and folding~\cite{Barsanti}).

\subsection{Preliminary remarks on a single filament\label{introtech}}

The rigorous definitions will be given in section 2. Here we introduce less
formally a few objects related to a \textit{single} vortex filament.

We consider a 3d-Brownian motion $\{X_{t}\}_{t\in\lbrack0,T]}$ starting from a
point $X_{0}$. This is the backbone of the vortex filament whose vorticity
field is given by
\begin{equation}
\xi_{\text{single}}(x)=\frac{U}{\ell^{2}}\int_{0}^{T}\rho_{\ell}(x-X_{t})\circ
dX_{t}.
\end{equation}
The letter $t$, that sometimes we shall also call time, is not physical time
but just the parameter of the curve. All our random fields are time
independent, in the spirit of equilibrium statistical mechanics. We assume
that $\rho_{\ell}(x)=\rho(x/\ell)$ for a radially symmetric measurable bounded
function $\rho$ with compact support in the ball $B(0,1)$ (the unit ball in 3d
Euclidean space). To have an idea, consider the case $\rho=1_{B(0,1)}$. Then
$\xi_{\text{single}}(x)=0$ outside an $\ell$-neighboor $\mathcal{U}_{\ell}$ of
the support of the curve $\{X_{t}\}_{t\in\lbrack0,T]}$. Inside $\mathcal{U}%
_{\ell}$, $\xi_{\text{single}}(x)$ is a time-average of the ``directions''
$dX_{t}$, with the pre-factor ${U}/{\ell^{2}}$. More precisely, if
$x\in\mathcal{U}_{\ell}$, one has to consider the time-set where $X_{t}\in
B(x,\ell)$ and average $dX_{t}$ on such time set. The resulting field
$\xi_{\text{single}}(x)$ looks much less irregular than $\{X_{t}%
\}_{t\in\lbrack0,T]}$, with increasing irregularity for smaller values of
$\ell$.

When $\rho$ is only measurable, it is not a priori clear that the Stratonovich
integral $\xi_{\text{single}}(x)$ is well defined, since the quadratic
variation corrector involves distributional derivatives of $\rho$ (the It\^{o}
integral is more easily defined, but it is not the natural object to be
considered, see the remarks in \cite{Fl1}). Since we shall never use
explicitly $\xi_{\text{single}}(x)$, this question is secondary and we may
consider $\xi_{\text{single}}(x)$ just as a formal expression that we
introduce to motivate the subsequent definition of $u_{\text{single}}(x)$.
However, at least in some particular case ($\rho=1_{B(0,1)}$) or under a
little additional assumptions, the Stratonovich integral $\xi_{\text{single}%
}(x)$ is well defined since the corrector has a meaning (in the case of
$\rho=1_{B(0,1)}$ the corrector involves the local time of 3d Brownian motion
on sferical surfaces).

The factor ${U}/{\ell^{2}}$ in the definition of $\xi_{\text{single}}(x)$ is
obscure at this level. Formally, it could be more natural just to introduce a
parameter $\Gamma$, in place of ${U}/{\ell^{2}}$, to describe the intensity of
the vortex. However, we do not have a clear interpretation of $\Gamma$, a
posteriori, from out theorems, while on the contrary it will arise that the
parameter $U$ has the meaning of a typical velocity intensity in the most
active region of the filament. Thus the choice of the expression ${U}%
/{\ell^{2}}$ has been devised a posteriori. The final interpretation of the
three parameters is that $T$ is the length of the filament, $\ell$ the
thickness, $U$ the typical velocity around the core.

%
%
%
%
%
%
%
%

The velocity field $u$ generated by $\xi$ is given by the Biot-Savart
relation
\begin{equation}
u_{\text{single}}(x)=\frac{U}{\ell^{2}}\int_{0}^{T}K_{\ell}(x-X_{t}%
)\wedge\circ dX_{t}%
\end{equation}
where the vector kernel $K_{\ell}(x)$ is defined as
\begin{equation}
K_{\ell}(x)=\nabla V_{\ell}(x)=\frac{1}{4\pi}\int_{B(0,\ell)}\rho_{\ell
}(y)\frac{x-y}{|x-y|^{3}}dy,\qquad V_{\ell}(x)=-\frac{1}{4\pi}\int_{B(0,\ell
)}\frac{\rho_{\ell}(y)dy}{|x-y|}.
\end{equation}

The scalar field $V_{\ell}$ satisfy the Poisson equation $\Delta V_{\ell}%
=\rho_{\ell}$ in all $\mathbb{R}^{3}$. Since $\rho$ is radially symmetric and
with compact support we can have two different situations according to the
fact that the integral of $\rho$: $Q=\int_{B(0,1)}\rho(x)dx$ is zero or not.
If $Q=0$ then the field $V_{\ell}$ is identically zero outside the ball
$B(0,\ell)$. Otherwise the fields $V_{\ell},K_{\ell}$ outside the ball
$B(0,\ell)$ have the form
\begin{equation}
V_{\ell}(x)=Q\frac{\ell^{3}}{|x|},\quad K_{\ell}(x)=-2Q\ell^{3}\frac
{x}{|x|^{3}},\qquad\text{$|x|\geq\ell$}. \label{VKscal}%
\end{equation}
Accordingly we will call the case $Q=0$ \emph{short range} and $Q\neq0$
\emph{long range}. The proof of the previous formula for $K$ is given, for
completeness, in the remark at the end of the section.

A basic result is that the Stratonovich integral in the expression of
$u_{\text{single}}$ can be replaced by an It\^{o} integral:

\begin{lemma}
\label{eq:strato-ito} It\^{o} and Stratonovich integrals in the definition of
$u_{\text{single}}(x)$ coincide:
\begin{equation}
u_{\text{single}}(x)=\frac{U}{\ell^{2}}\int_{0}^{T}K_{\ell}(x-X_{t})\wedge
dX_{t}\label{eq:single-vortex-u-ito}%
\end{equation}
where the integral is understood in It\^{o} sense. Then $u_{\text{single}}$ is
a (local)-martingale with respect to the standard filtration of $X$.
\end{lemma}

About the proof, by an approximation procedure that we omit we may assume
$\rho$ H\"{o}lder continuous, so the derivatives of $K_{\ell}$ exist and are
H\"{o}lder continuous by classical regularity theorems for elliptic equations.
Under this regularity one may compute the corrector and prove that it is equal
to zero, so, a posteriori, equation (\ref{eq:single-vortex-u-ito}) holds true
in the limit also for less regular $\rho$. About the proof that the corrector
is zero, it can be done component-wise, but it is more illuminating to write
the following heuristic computation:\ the corrector is formally given by%
\[
-\frac{1}{2}\int_{0}^{T}\left(  \nabla K_{\ell}(x-X_{t})dX_{t}\right)  \wedge
dX_{t}.
\]
Now, from the property $dX_{t}^{i}dX_{t}^{j}=\delta_{ij}dt$ one can verify
that
\[
\left(  \nabla K_{\ell}(x-X_{t})dX_{t}\right)  \wedge dX_{t}=\left(
\text{curl }K_{\ell}\right)  (x-X_{t})dt.
\]
Since $K_{\ell}$ is a gradient, we have curl $K_{\ell}=0$, so the corrector is
equal to zero.

\begin{remark}
One may verify that div$~u_{\text{single}}=0$, so $u_{\text{single}}$ may be
the velocity fluid of an incompressible fluid. On the contrary, div$~\xi
_{\text{single}}$ is different from zero and curl$~u_{\text{single}}$ is not
$\xi_{\text{single}}$ but its projection on divergence free fields. Therefore,
one should think of $\xi_{\text{single}}$\ as an auxiliary field we start from
in the construction of the model.
\end{remark}

\begin{remark}
\label{formulaK}Let us prove (\ref{VKscal}), limited to $K$ to avoid
repetitions. We want to solve $\Delta V=\rho$ with $\rho$ spherically
symmetric. The gradient $K=\nabla V$ of $V$ satisfies $\nabla\cdot K=\rho$
and, by spherical symmetry, it must be such that
\[
K(x)=\frac{x}{|x|}f(|x|)
\]
for some scalar function $f(r)$. By Gauss theorem we have
\[
\int_{B(0,r)}\nabla\cdot K(x)dx=\int_{\partial B(0,r)}K(x)\cdot d\sigma(x)
\]
where $d\sigma(x)$ is the outward surface element of the sphere. So
\[
f(r)4\pi r^{2}=\int_{B(0,r)}\rho(x)dx=:Q(r)
\]
namely $f(r)=\frac{Q(r)}{4\pi r^{2}}$. Therefore
\[
K(x)=Q(r)\frac{x}{4\pi|x|^{3}}%
\]
with $Q(r)=Q(1)$ if $r\geq1$ (since $\rho$ has support in $B(0,1)$).
\end{remark}

\section{Poisson field of vortices}

Intuitively, we want to describe a collection of infinitely many independent
Brownian vortex filaments, uniformly distributed in space, with
intensity-thickness-length parameters $(U,\ell,T)$ distributed according to a
measure $\gamma$. The total vorticity of the fluid is the sum of the vorticity
of the single filaments, so, by linearity of the relation vorticity-velocity,
the total velocity field will be the sum of the velocity fields of the single filaments.

The rigorous description requires some care, so we split it into a number of steps.

\subsection{Underlying Poisson random measure}

Let $\Xi$ be the metric space
\[
\Xi=\{(U,\ell,T,X)\in\mathbb{R}_{+}^{3}\times C([0,1];\mathbb{R}^{3}%
):0<\ell\leq\sqrt{T}\leq1\}
\]
with its Borel $\sigma$-field $\mathcal{B}\left(  \Xi\right)  $. Let $\left(
\Omega,\mathcal{A},P\right)  $ be a probability space, with expectation
denoted by $E$, and let $\mu_{\omega}$, $\omega\in\Omega$, be a Poisson random
measure on $\mathcal{B}\left(  \Xi\right)  $, with intensity $\nu$ (a $\sigma
$-finite measure on $\mathcal{B}\left(  \Xi\right)  $) given by
\[
d\nu(U,\ell,T,X)=d\gamma(U,\ell,T)d\mathcal{W}(X).
\]
for $\gamma$ a $\sigma$-finite measure on the Borel sets of $\{(U,\ell
,T)\in\mathbb{R}_{+}^{3}:0<\ell\leq\sqrt{T}\leq1\}$ and $d\mathcal{W}(X)$ the
$\sigma$-finite measure defined by
\[
\int_{C([0,1];\mathbb{R}^{3})} \psi(X) d\mathcal{W}(X) = \int_{\mathbb{R}^{3}}
\left[  \int_{C([0,1];\mathbb{R}^{3})} \psi(X) d\mathcal{W}_{x_{0}}(X)\right]
dx_{0}%
\]
for any integrable test function $\psi: C([0,1];\mathbb{R}^{3}) \to\mathbb{R}$
where $d\mathcal{W}_{x_{0}}(X)$ is the Wiener measure on $C([0,1],\mathbb{R}%
^{3})$ starting at $x_{0}$ and $dx_{0}$ is the Lebesgue measure on
$\mathbb{R}^{3}$. Heuristically the measure $\mathcal{W}$ describe a Brownian
path starting from an uniformly distributed point in all space. The
assumptions on $\gamma$ will be specified at due time.

The random measure $\mu_{\omega}$ is uniquely determined by its characteristic
function
\[
\expect\exp\left(  i\int_{\Xi}\varphi(\xi)\mu(d\xi)\right)  =\exp\left(
-\int_{\Xi}(e^{\varphi(\xi)}-1)\nu(d\xi)\right)
\]
for any bounded measurable function $\varphi$ on $\Xi$ with support in a set
of finite $\nu$-measure.

In particular, for example, the first two moments of $\mu$ read
\[
\expect\int_{\Xi}\varphi(\xi)\mu(d\xi)=\int_{\Xi}\varphi(\xi)\nu(d\xi)
\]
and
\[
\expect\left[  \int_{\Xi}\varphi(\xi)\mu(d\xi)\right]  ^{2}=\left[  \int_{\Xi
}\varphi(\xi)\nu(d\xi)\right]  ^{2}+\int_{\Xi}\varphi^{2}(\xi)\nu(d\xi).
\]
We have to deal also with moments of order $p$;\ some useful formulae will be
now given.

\subsection{Moments of the Poisson Random Field}

Let $\varphi:\Xi\rightarrow\mathbb{R}$ be a measurable function. We shall say
it is $\mu$-integrable if it is $\mu_{\omega}$-integrable for $P$-a.e.
$\omega\in\Omega$. In such a case, by approximation by bounded measurable
compact support functions, one can show that the mapping $\omega\mapsto
\mu_{\omega}\left(  \varphi\right)  $ is measurable.

If $\varphi^{\prime}:\Xi\rightarrow\mathbb{R}$ is a measurable function with
$\nu\left(  \varphi\neq\varphi^{\prime}\right)  =0$, since $P\left(
\mu\left(  \varphi\neq\varphi^{\prime}\right)  =0\right)  =1$, we have that
$\varphi^{\prime}$ is $\mu$-integrable and $P\left(  \mu\left(  \varphi
\right)  =\mu\left(  \varphi^{\prime}\right)  \right)  =1$. Therefore the
concept of $\mu$-integrability and the random variable $\mu\left(
\varphi\right)  $ depend only on the equivalence class of $\varphi$.

Let $\varphi$ be a measurable function on $\Xi$, possibly defined only $\nu
$-a.s. We shall say that it is $\mu$-integrable if some of its measurable
extensions to the whole $\Xi$ is $\mu$-integrable. Neither the condition of
$\mu$-integrability nor the equivalence class of $\mu\left(  \varphi\right)  $
depend on the extension, by the previous observations. We need all these
remarks in the sequel when we deal with $\varphi$ given by stochastic integrals.

\begin{lemma}
\label{moments}Let $\varphi$ be a measurable function on $\Xi$ (possibly
defined only $\nu$-a.s.), such that $\nu(\varphi^{p})<\infty$ for some even
integer number $p$. Then $\varphi$ is $\mu$-integrable, $\mu(\varphi)\in
L^{p}\left(  \Omega\right)  $, and
\[
E\left[  \mu\left(  \varphi\right)  ^{p}\right]  \leq ep^{p}\nu(\varphi^{p}).
\]
If in addition we have $\nu(\varphi^{k})=0$ for every odd $k<p$, then
\[
E\left[  \mu\left(  \varphi\right)  ^{p}\right]  \geq\nu(\varphi^{p}).
\]
\end{lemma}

\textsc{Proof. }Assume for a moment that $\varphi$ is bounded measurable and
with support in a set of finite $\nu$-measure, so that the qualitative parts
of the statement are obviously true. Using the moment generating function
\[
\expect e^{\lambda\mu(\varphi)}=e^{\nu(e^{\lambda\varphi}-1)}%
\]
we obtain
\[%
\begin{split}
\sum_{p=0}^{\infty}\frac{\lambda^{p}}{p!}\expect[ \mu(\varphi)^p ]  &
=\sum_{n=0}^{\infty}\frac{1}{n!}\left[  \nu(e^{\lambda\varphi}-1)\right]
^{n}\\
&  =1+\sum_{n=1}^{\infty}\frac{1}{n!}\left[  \sum_{k=1}^{\infty}\frac
{\lambda^{k}}{k!}\nu(\varphi^{k})\right]  ^{n}\\
&  =1+\sum_{n=1}^{\infty}\frac{1}{n!}\sum_{k_{1},\dots,k_{n}\geq1}%
\frac{\lambda^{k_{1}+\dots+k_{n}}}{k_{1}!\cdots k_{n}!}\nu(\varphi^{k_{1}%
})\cdots\nu(\varphi^{k_{n}})\\
&  =1+\sum_{p=1}^{\infty}\frac{\lambda^{p}}{p!}\sum_{n=1}^{p}\sum
_{\substack{k_{1},\dots,k_{n}\geq1 \\k_{1}+\cdots+k_{n}=p}}\frac{p!}%
{n!k_{1}!\cdots k_{n}!}\nu(\varphi^{k_{1}})\cdots\nu(\varphi^{k_{n}})
\end{split}
\]
Hence we have an equation for the moments:
\begin{equation}
\expect[ \mu(\varphi)^p ]=\sum_{n=1}^{p}\sum_{\substack{k_{1},\dots,k_{n}\geq1
\\k_{1}+\cdots+k_{n}=p}}\frac{p!}{n!k_{1}!\cdots k_{n}!}\nu(\varphi^{k_{1}%
})\cdots\nu(\varphi^{k_{n}}) \label{sum}%
\end{equation}
Since $\nu(|\varphi|^{k})\leq\lbrack\nu(|\varphi|^{p})]^{k/p}$ for $k\leq p$
we have
\[%
\begin{split}
\expect[|\mu(\varphi)|^p]  &  \leq\expect[ \mu(|\varphi|)^p ]\\
&  \leq\sum_{n=1}^{p}\sum_{\substack{k_{1},\dots,k_{n}\geq1 \\k_{1}%
+\cdots+k_{n}=p}}\frac{p!}{n!k_{1}!\cdots k_{n}!}\nu(|\varphi|^{p})\\
&  \leq\nu(|\varphi|^{p})\sum_{n=1}^{p}\sum_{\substack{k_{1},\dots,k_{n}\geq0
\\k_{1}+\cdots+k_{n}=p}}\frac{p!}{n!k_{1}!\cdots k_{n}!}\\
&  =\nu(|\varphi|^{p})\sum_{n=1}^{p}\frac{n^{p}}{n!}\leq ep^{p}\nu
(|\varphi|^{p}).
\end{split}
\]
This proves the first inequality of the lemma. A posteriori, we may use it to
prove the qualitative part of the first statement, by a simple approximation
procedure for general measurable $\varphi$.

For the lower bound, from the assumption that $\nu(\varphi^{k})=0$ for odd
$k<p$, in the sum (\ref{sum}) we have contributions only when all $k_{i}$,
$i=1,\dots,n$ are even. Then, neglecting many terms, we have
\[
\expect[ \mu(\varphi)^p ]\geq\nu(\varphi^{p}).
\]
The proof is complete.

\subsection{Velocity field}

Let $\rho$ be a radially symmetric measurable bounded function $\rho$ on
$\mathbb{R}^{3}$ with compact support in the ball $B(0,1)$, and let $K_{\ell}$
be defined as in section \ref{introtech}. Then we have that
\begin{equation}
K_{1}\text{ is Lipschitz continuous. } \label{lip}%
\end{equation}
Indeed by an explicit computation it is possible to show that
\[
K_{1}(x) = - 2 Q(|x|) \frac{x}{|x|^{3}}
\]
with
\[
Q(r) := \int_{B(0,r)} \rho(x) dx
\]
(since $\rho$ has support in $B(0,1)$ we have $Q(r) = Q$ if $r \ge1$). Then
\[
\nabla K_{1}(x) = - 2 Q^{\prime}(|x|) \frac{x \otimes x}{|x|^{4}} - 2 Q(|x|)
\left[  \frac{1}{|x|^{3}}-3\frac{x\otimes x}{|x|^{5}}\right]
\]
with $Q^{\prime}(r) = dQ(r)/dr$ and we can bound
\[
|\nabla K_{1}(x)| \le C \|\rho\|_{\infty}%
\]
since
\[
Q(r) \le C \|\rho\|_{\infty}r^{3}, \qquad Q^{\prime}(r) \le C \|\rho
\|_{\infty}r^{2}.
\]

\medskip
%
%
%

For any $x,x_{0}\in\mathbb{R}^{3}$, $\ell,T>0$, the random variable
$X\mapsto\int_{0}^{T}K_{\ell}(x-X_{t})\wedge dX_{t}$ is defined $\mathcal{W}%
_{x_{0}}$-a.s. on $C([0,1],\mathbb{R}^{3})$. We also have, given
$x\in\mathbb{R}^{3}$, $\ell,T>0$, that $X\mapsto\int_{0}^{T}K_{\ell}%
(x-X_{t})\wedge dX_{t}$ is a well defined measurable function, defined
$\mathcal{W}$-a.s. on $C([0,1],\mathbb{R}^{3})$. More globally, writing
$\xi=(U,\ell,T,X)$ for shortness, for any $x\in\mathbb{R}^{3}$ we may consider
the measurable $\mathbb{R}^{3}$-valued function
\[
\xi\mapsto u_{\text{single}}^{\xi}\left(  x\right)  :=\frac{U}{\ell^{2}}%
\int_{0}^{T}K_{\ell}(x-X_{t})\wedge dX_{t}%
\]
defined $\mathcal{\nu}$-a.s. on $\Xi$. In plain words, this is the velocity
field at point $x$ of a filament specified by $\xi$.

Again in plain words, given $\omega\in\Omega$, the point measure $\mu_{\omega
}$ specifies the parameters and locations of infinitely many filaments:
formally
\begin{equation}
\mu=\sum_{\alpha\in\mathbb{N}}\delta_{\xi^{\alpha}} \label{muformal}%
\end{equation}
for a sequence of i.i.d. random points $\left\{  \xi^{\alpha}\right\}  $
distributed in $\Xi$ according to $\nu$ (this fact is not rigorous since $\nu$
is only $\sigma$-finite, but it has a rigorous version by localization
explained below). Since the total velocity at a given point $x\in
\mathbb{R}^{3}$ should be the sum of the contributions from each single
filament, i.e. in heuristic terms
\begin{equation}
u(x)=\sum_{\alpha\in\mathbb{N}}u_{\text{single}}^{\xi^{\alpha}}\left(
x\right)  \label{uformal}%
\end{equation}
we see that, in the rigorous language of $\mu$, we should write
\begin{equation}
u(x)=\int_{\Xi}u_{\text{single}}^{\xi}\left(  x\right)  \mu(d\xi)=\mu\left(
u_{\text{single}}^{\cdot}\left(  x\right)  \right)  .
\end{equation}
If we show that $\nu\left(  \left|  u_{\text{single}}^{\cdot}\left(  x\right)
\right|  ^{p}\right)  <\infty$ for some even $p\geq2$, then from lemma
\ref{moments}, $u_{\text{single}}^{\cdot}\left(  x\right)  $ is $\mu
$-integrable and the random variable
\[
\omega\mapsto u(x,\omega):=\mu_{\omega}\left(  u_{\text{single}}^{\cdot
}\left(  x\right)  \right)
\]
is well defined.

Since the proof of the next lemma we start to use the occupation measure of
the 3d Brownian motion. We set
\[
L_{B}^{T}:=\int_{0}^{T}1_{X_{t}\in B}dt
\]
for every Borel set $B$ of $\mathbb{R}^{3}$.

\begin{lemma}
\label{singlebound}Given $x\in\mathbb{R}^{3}$ and $p>0$, there exist $C_{p}>0$
such that for every $\ell^{2}\leq T\leq1$we have
\[
\mathcal{W}\left[  W(x)^{p/2}\right]  \leq C_{p}\ell^{2p}\ell T
\]
where
\[
W(x)=\int_{0}^{T}\left|  K_{\ell}(x-X_{t})\right|  ^{2}dt.
\]
\end{lemma}

\textsc{Proof. }Let us bound $K$ by a multi-scale argument. This is necessary
only in the long-range case (see the introduction). If $|y|\leq\ell$ we can
bound $|K_{\ell}(y)|\leq C\ell$. Indeed if $|w|\leq1$, $|K(w)|\leq C$ for some
constant $C$ and then if $|y|<\ell$ we have $|K_{\ell}(y)|=\ell|K_{1}%
(y/\ell)|\leq C\ell$.

Next, given $\Lambda>\ell$ and an integer $N$, consider a sequence $\left\{
\ell_{i}\right\}  _{i=0,...,N}$ of scales, with $\ell=\ell_{0}<\ell
_{1}<...<\ell_{N}=\Lambda$. Then, for $i=1,\dots,N$, if $\ell_{i-1}%
<|y|<\ell_{i}$, by the explicit formula for $K_{\ell}(y)$ we have $|K_{\ell
}(y)|\leq C\ell(\ell/|y|)^{2}$. Therefore
\[
|K_{\ell}(y)|\leq C\ell\left(  \frac{\ell}{\ell_{i-1}}\right)  ^{2}%
1_{\ell_{i-1}<|y|\leq\ell_{i}}%
\]

If $|y|>\Lambda$ we simply bound $K_{\ell}^{e}(y)\leq C\ell(\ell/\Lambda)^{2}
$.

Summing up,
\begin{align*}
\left|  K_{\ell}(y)\right|  ^{2}  &  =\left|  K_{\ell}(y)\right|  ^{2}\left(
1_{|y|\leq\ell}+\sum_{i=1}^{N}1_{\ell_{i-1}<|y|<\ell_{i}}+1_{|y|>\Lambda
}\right) \\
&  \leq C\ell^{2}1_{|y|\leq\ell}+C\sum_{i=1}^{N}\ell^{2}\left(  \frac{\ell
}{\ell_{i-1}}\right)  ^{4}1_{\ell_{i-1}<|y|<\ell_{i}}+C\ell^{2}\left(
\frac{\ell}{\Lambda}\right)  ^{4}1_{|y|>\Lambda}%
\end{align*}
which implies the following bound for $W(x)$:
\[
W(x)\leq C\ell^{2}L_{B(x,\ell)}^{T}+C\sum_{i=1}^{N}\ell^{2}\left(  \frac{\ell
}{\ell_{i-1}}\right)  ^{4}L_{B(x,\ell_{i})\backslash B(x,\ell_{i-1})}%
^{T}+C\ell^{2}\left(  \frac{\ell}{\Lambda}\right)  ^{4}L_{B(x,\Lambda)^{c}%
}^{T}%
\]
where $L_{B}^{T}$ has been defined above. By the additivity of $B\mapsto
L_{B}^{T}$, the sum appearing in this equation can be rewritten as
\begin{align*}
\sum_{i=1}^{N}\left(  \frac{\ell}{\ell_{i-1}}\right)  ^{4}L_{B(x,\ell
_{i})\backslash B(x,\ell_{i-1})}^{T}  &  =\sum_{i=1}^{N-1}\left[  \left(
\frac{\ell}{\ell_{i-1}}\right)  ^{4}-\left(  \frac{\ell}{\ell_{i}}\right)
^{4}\right]  L_{B(x,\ell_{i})}^{T}\\
&  +\left(  \frac{\ell}{\ell_{N-1}}\right)  ^{4}L_{B(x,\Lambda)}%
^{T}-L_{B(x,\ell)}^{T}%
\end{align*}%
\[
\leq\sum_{i=1}^{N-1}\left[  \left(  \frac{\ell}{\ell_{i-1}}\right)
^{4}-\left(  \frac{\ell}{\ell_{i}}\right)  ^{4}\right]  L_{B(x,\ell_{i})}%
^{T}+\left(  \frac{\ell}{\ell_{N-1}}\right)  ^{4}L_{B(x,\Lambda)}^{T}.
\]
Assume that $\ell_{i}/\ell_{i-1}\leq2$ uniformly in $i=1,\dots,N$. Then
\begin{align*}
&  \sum_{i=1}^{N}\left(  \frac{\ell}{\ell_{i-1}}\right)  ^{4}L_{B(x,\ell
_{i})\backslash B(x,\ell_{i-1})}^{T}\\
&  \leq C\sum_{i=1}^{N-1}\left[  \left(  \frac{\ell}{\ell_{i-1}}\right)
^{2}-\left(  \frac{\ell}{\ell_{i}}\right)  ^{2}\right]  \left(  \frac{\ell
}{\ell_{i}}\right)  ^{2}L_{B(x,\ell_{i})}^{T}+C\left(  \frac{\ell}{\Lambda
}\right)  ^{4}L_{B(x,\Lambda)}^{T}.
\end{align*}

Notice now that
\[
\sum_{i=1}^{N-1}\left[  \left(  \frac{\ell}{\ell_{i-1}}\right)  ^{2}-\left(
\frac{\ell}{\ell_{i}}\right)  ^{2}\right]  =\left(  \frac{\ell}{\ell}\right)
^{2}-\left(  \frac{\ell}{\Lambda}\right)  ^{2}\leq1
\]
so that, by Cauchy-Schwartz and Jensen inequalities we have
\begin{align*}
&  \left[  \sum_{i=1}^{N}\left(  \frac{\ell}{\ell_{i-1}}\right)
^{4}L_{B(x,\ell_{i})\backslash B(x,\ell_{i-1})}^{T}\right]  ^{p/2}\\
&  \leq C_{p}\sum_{i=1}^{N-1}\left[  \left(  \frac{\ell}{\ell_{i-1}}\right)
^{2}-\left(  \frac{\ell}{\ell_{i}}\right)  ^{2}\right]  \left(  \frac{\ell
}{\ell_{i}}\right)  ^{p}(L_{B(x,\ell_{i})}^{T})^{p/2}+C_{p}\left(  \frac{\ell
}{\Lambda}\right)  ^{2p}\left(  L_{B(x,\Lambda)}^{T}\right)  ^{p/2}.
\end{align*}
An upper bound for $W(x)^{p/2}$ is then obtained as
\begin{align*}
W(x)^{p/2}  &  \leq C_{p}\ell^{p}(L_{B(x,\ell)}^{T})^{p/2}\\
&  +C_{p}\ell^{p}\left(  \frac{\ell}{\Lambda}\right)  ^{2p}\left(
L_{B(x,\Lambda)}^{T}\right)  ^{p/2}+C_{p}\ell^{p}\left(  \frac{\ell}{\Lambda
}\right)  ^{2p}\left(  L_{B(x,\Lambda)^{c}}^{T}\right)  ^{p/2}\\
&  +C\ell^{p}\sum_{i=1}^{N-1}\left[  \left(  \frac{\ell}{\ell_{i-1}}\right)
^{2}-\left(  \frac{\ell}{\ell_{i}}\right)  ^{2}\right]  \left(  \frac{\ell
}{\ell_{i}}\right)  ^{p}(L_{B(x,\ell_{i})}^{T})^{p/2}%
\end{align*}%
\begin{align*}
&  \leq C_{p}\ell^{p}(L_{B(x,\ell)}^{T})^{p/2}+C_{p}\ell^{p}\left(  \frac
{\ell}{\Lambda}\right)  ^{2p}T^{p/2}\\
&  +C_{p}\ell^{p}\sum_{i=1}^{N-1}\left[  \left(  \frac{\ell}{\ell_{i-1}%
}\right)  ^{2}-\left(  \frac{\ell}{\ell_{i}}\right)  ^{2}\right]  \left(
\frac{\ell}{\ell_{i}}\right)  ^{p}(L_{B(x,\ell_{i})}^{T})^{p/2}%
\end{align*}
where we have used again Cauchy-Schwartz inequality.

We use now lemma~\ref{lemma:definitive-L-bounds} with $\alpha=1$. For a given
$\lambda\in\left(  0,1\right)  $, we take bot $\varepsilon$ and $\ell$ equal
to $\lambda$ in (\ref{eq:Lmoments-bound4}) and (\ref{eq:Lmoments-bound4-bis}),
and get
\[
\mathcal{W}\left[  \left(  L_{B(x,\lambda)}^{T}\right)  ^{p/2}\right]  \leq
C_{p}(\lambda\wedge\sqrt{T})^{p}\lambda(\lambda\vee\sqrt{T})^{2}.
\]
Then we obtain
\begin{align*}
&  \mathcal{W}\left[  W(x)^{p/2}\right] \\
&  \leq C_{p}\ell^{2p}\ell T+C_{p}\ell^{p}\left(  \frac{\ell}{\Lambda}\right)
^{2p}T^{p/2}\\
&  +C_{p}\ell^{p}\sum_{i=1}^{N-1}\left[  \left(  \frac{\ell}{\ell_{i-1}%
}\right)  ^{2}-\left(  \frac{\ell}{\ell_{i}}\right)  ^{2}\right]  \left(
\frac{\ell}{\ell_{i}}\right)  ^{p}(\ell_{i}\wedge\sqrt{T})^{p}\ell_{i}%
(\ell_{i}\vee\sqrt{T})^{2}%
\end{align*}
and taking the limit as the partition gets finer:
\[%
\begin{split}
&  \mathcal{W}\left[  W(x)^{p/2}\right] \\
&  \leq C_{p}\ell^{2p}\ell T+C_{p}\ell^{p}\left(  \frac{\ell}{\Lambda}\right)
^{2p}T^{p/2}\\
&  +C_{p}\ell^{p}\int_{\ell}^{\Lambda}\left(  \frac{\ell}{u}\right)
^{p}(u\wedge\sqrt{T})^{p}(u\vee\sqrt{T})^{2}ud\left[  -\left(  \frac{\ell}%
{u}\right)  ^{2}\right]
\end{split}
\]
The integral can then be computed as
\[%
\begin{split}
\int_{\ell}^{\Lambda}  &  \left(  \frac{\ell}{u}\right)  ^{p}(u\wedge\sqrt
{T})^{p}(u\vee\sqrt{T})^{2}u\frac{\ell^{2}}{u^{3}}du\\
&  =\ell^{p}\ell^{2}T\int_{\ell}^{\sqrt{T}}\frac{du}{u^{2}}+T^{p/2}\ell
^{p+2}\int_{\sqrt{T}}^{\Lambda}u^{-p}du\\
&  =\ell^{p}\ell^{2}T[\ell^{-1}-\sqrt{T}^{-1}]+(p-1)^{-1}T^{p/2}\ell
^{p+2}[T^{(1-p)/2}-\Lambda^{(1-p)}]\\
&  \leq\ell^{p}\ell T+(p-1)^{-1}\ell^{p}\frac{\ell}{\sqrt{T}}\ell T
\end{split}
\]
Using the fact that $\ell\leq\sqrt{T}$ and letting $\Lambda\rightarrow\infty$
we finally obtain the claim. The proof is complete. \hfill$\square$\\[0.5cm]

\begin{remark}
The multiscale argument above can be rewritten in contunuum variables from the
very beggining by means of the following identity: if $f:[0,\infty
)\rightarrow\mathbb{R}$ is of class $C^{1}$ and has a suitable decay at
infinity, then%
\[
\int_{0}^{T}f\left(  \left|  x-X_{t}\right|  \right)  dt=-\int_{0}^{\infty
}f^{\prime}\left(  r\right)  L_{B(x,r)}^{T}dr.
\]
This identity can be applied to $W(x)$. The proof along these lines is not
essentially shorter and perhaps it is more obscure, thus we have choosen the
discrete multiscale argument which has a neat geometrical interpretation.
\end{remark}

\begin{corollary}
\label{corollintegr}Assume
\[
\gamma(U^{p}\ell T)<\infty
\]
for some even integer $p\geq2$. Then, for any $x\in\mathbb{R}^{3}$, we have
\[
\nu\left(  \left|  u_{\text{single}}^{\cdot}\left(  x\right)  \right|
^{p}\right)  <\infty
\]
and the random variable
\[
\omega\mapsto u(x,\omega):=\mu_{\omega}\left(  u_{\text{single}}^{\cdot
}\left(  x\right)  \right)
\]
has finite $p$-moment:
\[
E\left[  \left|  u(x)\right|  ^{p}\right]  <\infty.
\]
\end{corollary}

\textsc{Proof. }We have
\[
\nu\left(  \left|  u_{\text{single}}^{\cdot}\left(  x\right)  \right|
^{p}\right)  =\int\mathcal{W}_{x_{0}}(\left|  u_{\text{single}}^{\cdot}\left(
x\right)  \right|  ^{p})d\gamma(U,\ell,T)dx_{0}.
\]
By Burkholder-Davis-Gundy inequality, there is $C_{p}>0$ such that
\[
\mathcal{W}_{x_{0}}(\left|  u_{\text{single}}^{\cdot}\left(  x\right)
\right|  ^{p})\leq C_{p}\frac{U^{p}}{\ell^{2p}}\mathcal{W}_{x_{0}}\left[
W(x)^{p/2}\right]  .
\]
Hence
\begin{align*}
\mathcal{W}(\left|  u_{\text{single}}^{\cdot}\left(  x\right)  \right|  ^{p})
&  \leq C_{p}\frac{U^{p}}{\ell^{2p}}\mathcal{W}\left[  W(x)^{p/2}\right] \\
&  \leq C_{p}^{\prime}\frac{U^{p}}{\ell^{2p}}\ell^{2p}\ell T=C_{p}^{\prime
}U^{p}\ell T.
\end{align*}
Therefore $\nu\left(  \left|  u_{\text{single}}^{\cdot}\left(  x\right)
\right|  ^{p}\right)  <\infty$ by the assumption $\gamma(U^{p}\ell T)<\infty$.
The other claims are a consequence of lemma \ref{moments}. \hfill$\square$

\begin{lemma}
\label{symmetries}Under the previous assumptions, the law of $u(x,\cdot)$ is
independent of $x$ and is invariant also under rotations:
\[
Ru(x,\cdot)\overset{\mathcal{L}}{=}u(Rx,\cdot)
\]
for every rotation matrix $R$.
\end{lemma}

\textsc{Proof. }With the usual notation $\xi=(U,\ell,T,X)$ we have
\[
u_{\text{single}}^{\xi}\left(  x\right)  =u_{\text{single}}^{(U,\ell
,T,X)}\left(  x\right)  =u_{\text{single}}^{(U,\ell,T,X-x)}\left(  0\right)
=u_{\text{single}}^{\tau_{x}\xi}\left(  0\right)
\]
where $\tau_{x}(U,\ell,T,X)=(U,\ell,T,X-x)$. The map $\tau_{x}$ is a
measurable transformation of $\Xi$ into itself. One can see that $\nu$ is
$\tau_{x}$-invariant; we omit the details, but we just notice that $\nu$ is
not a finite measure, so the invariance means
\[
\int_{\Xi}\varphi\left(  \tau_{x}\xi\right)  \nu\left(  d\xi\right)
=\int_{\Xi}\varphi\left(  \xi\right)  \nu\left(  d\xi\right)
\]
for every $\varphi\in L^{1}\left(  \Xi,\nu\right)  $. From this invariance it
follows that the law of the random measure $\mu$ is the same as the law of the
random measure $\tau_{x}\mu$. Therefore
\begin{align*}
\mu\left[  u_{\text{single}}^{\xi}\left(  x\right)  \right]   &  =\mu\left[
u_{\text{single}}^{\tau_{x}\xi}\left(  0\right)  \right] \\
&  =\left(  \tau_{x}\mu\right)  \left[  u_{\text{single}}^{\xi}\left(
0\right)  \right]  \overset{\mathcal{L}}{=}\mu\left[  u_{\text{single}}^{\xi
}\left(  0\right)  \right]  .
\end{align*}
This proves the first claim.

If $R$ is a rotation, from the explicit form of $K_{\ell}$ it is easy to see
that
\[
RK_{\ell}(y)=K_{\ell}(Ry)
\]
hence
\begin{align*}
Ru_{\text{single}}^{\xi}\left(  x\right)   &  =\frac{U}{\ell^{2}}\int_{0}%
^{T}RK_{\ell}(x-X_{t})\wedge dRX_{t}\\
&  =\frac{U}{\ell^{2}}\int_{0}^{T}K_{\ell}(Rx-RX_{t})\wedge dRX_{t}\\
&  =u_{\text{single}}^{R\xi}\left(  Rx\right)
\end{align*}
where we have set $R(U,\ell,T,X)=(U,\ell,T,RX)$. Again $R\nu=\nu$,
$R\mu\overset{\mathcal{L}}{=}\mu$, so the end of the proof is the same as
above. \hfill$\square$

We say that a random field $u(x,\cdot)$ is homogeneous if it is independent of
$x$ and it is isotropic if it is invariant under rotations.

\begin{corollary}
\label{corollarybase}Assume
\[
\gamma(U^{p}\ell T)<\infty
\]
for every $p>1$. Then $\left\{  u(x,\cdot);x\in\mathbb{R}^{3}\right\}  $ is an
isotropic homogeneous random field, with finite moments of all orders.
\end{corollary}

This corollary is sufficient to introduce the structure function and state the
main results of this paper. However, it is natural to ask whether the random
field $\left\{  u(x,\cdot);x\in\mathbb{R}^{3}\right\}  $ has a continuous
modification. Having in mind Kolmogorov regularity theorem, we need good
estimates of $E\left[  \left|  u(x)-u(y)\right|  ^{p}\right]  $. They are as
difficult as the careful estimates we shall perform in the next section to
understand the scaling of the structure function. Therefore we anticipate the
result without proof. It is a direct consequence of Theorem \ref{maintheorem}.

\begin{proposition}
Assume
\[
\gamma(U^{p}\ell T)<\infty
\]
for every $p>1$. Then, for every even integer $p$ there is a constant
$C_{p}>0$ such that
\[
E\left[  \left|  u(x)-u(y)\right|  ^{p}\right]  \leq C_{p}\gamma\left[
U^{p}\left(  \frac{\ell\wedge\left|  x-y\right|  }{\ell}\right)  ^{p}\ell
T\right]  .
\]
Consequently, if the measure $\gamma$ has the property that for some even
integer $p$ and real number $\alpha>3$ there is a constant $C_{p}^{\prime}>0$
such that
\begin{equation}
\gamma\left[  U^{p}\left(  \frac{\ell\wedge\varepsilon}{\ell}\right)  ^{p}\ell
T\right]  \leq C_{p}^{\prime}\varepsilon^{\alpha}\text{ for any }%
\varepsilon\in\left(  0,1\right)  , \label{contgamma}%
\end{equation}
then the random field $u(x)$ has a continuous modification.
\end{proposition}

\begin{remark}
\label{remregularity}A sufficient condition for (\ref{contgamma}) is: there
are $\alpha>3$ and $\beta>0$ such that for every sufficiently large even
integer $p$ there is a constant $C_{p}>0$ such that
\[
\gamma\left[  U^{p}\ell T\cdot1_{\ell\leq\varepsilon^{1-\beta}}\right]  \leq
C_{p}\varepsilon^{\alpha}\text{ for any }\varepsilon\in\left(  0,1\right)  .
\]
Indeed,
\begin{align*}
&  \gamma\left[  U^{p}\left(  \frac{\ell\wedge\varepsilon}{\ell}\right)
^{p}\ell T\right] \\
&  =\gamma\left[  U^{p}\left(  \frac{\ell\wedge\varepsilon}{\ell}\right)
^{p}\ell T\cdot1_{\ell\leq\varepsilon^{1-\beta}}\right]  +\gamma\left[
U^{p}\left(  \frac{\ell\wedge\varepsilon}{\ell}\right)  ^{p}\ell
T\cdot1_{\varepsilon^{1-\beta}\leq\ell}\right] \\
&  \leq\gamma\left[  U^{p}\ell T\cdot1_{\ell\leq\varepsilon^{1-\beta}}\right]
+\gamma\left[  U^{p}\varepsilon^{\beta p}\ell T\cdot1_{\varepsilon^{1-\beta
}\leq\ell}\right] \\
&  \leq C_{p}\varepsilon^{\alpha}+\varepsilon^{\beta p}\gamma\left[  U^{p}\ell
T\right]
\end{align*}
so we have (\ref{contgamma}) with a suitable choice of $p$. This happens in
particular in the multifractal example of section \ref{multifrexample}, remark
\ref{remmultmodif}.
\end{remark}

The model presented here has a further symmetry which is not physically
correct. This symmetry, described in the next lemma, implies that the odd
moments of the longitudinal structure function vanish, in contradiction both
with experiments and certain rigorous results derived from the Navier-Stokes
equation (see \cite{frisch}). The same drawback is present in other
statistical models of vortex structures \cite{Pullin}.

Beyond the rigorous formulation, the following property says that the random
field $u_{\text{single}}^{\cdot}$ has the same law as $-u_{\text{single}%
}^{\cdot}$. We cannot use the concept of law since $\xi$ does not live on a
probability space.

\begin{lemma}
Given $x_{1},...,x_{n}\in\mathbb{R}^{3}$, the measurable vector
\[
U_{n}\left(  \xi\right)  :=\left(  u_{\text{single}}^{\xi}\left(
x_{1}\right)  ,...,u_{\text{single}}^{\xi}\left(  x_{n}\right)  \right)
\]
has the property
\[
\int\varphi\left(  U_{n}\left(  \xi\right)  \right)  d\nu\left(  \xi\right)
=\int\varphi\left(  -U_{n}\left(  \xi\right)  \right)  d\nu\left(
\xi\right)
\]
for every $\varphi=\varphi\left(  u_{1},...,u_{n}\right)  :\mathbb{R}%
^{3n}\rightarrow\mathbb{R}$ with a polynomial bound in its variables..
\end{lemma}

\textsc{Proof.}\hspace{0.3cm} With the notation $\widetilde{X}_{t}=X_{T-t}$,
we have
\begin{align*}
-u_{\text{single}}^{\xi}\left(  x_{k}\right)   &  =-\frac{U}{\ell^{2}}\int
_{0}^{T}K_{\ell}(x_{k}-X_{t})\wedge\circ dX_{t}\\
&  =\frac{U}{\ell^{2}}\int_{0}^{T}K_{\ell}(x_{k}-\widetilde{X}_{s})\wedge\circ
d\widetilde{X}_{s}\\
&  =u_{\text{single}}^{S\xi}\left(  x_{k}\right)
\end{align*}
where $S(U,\ell,T,X)=(U,\ell,T,\widetilde{X})$. Since $S\nu=\nu$, we have the
result (using the integrability of corollary \ref{corollintegr}).
\hfill$\square$

\begin{lemma}
\label{reflectsymm}If $p$ is an odd positive integer, then
\[
\nu\left[  \left\langle u_{\text{single}}^{\cdot}\left(  y\right)
-u_{\text{single}}^{\cdot}\left(  x\right)  ,y-x\right\rangle ^{p}\right]  =0
\]
for every $x,y\in\mathbb{R}^{3}$.
\end{lemma}

\textsc{Proof.}\hspace{0.3cm} It is sufficient to apply the lemma to the
function
\[
\varphi\left(  u_{1},u_{2}\right)  :=\left\langle u_{1}-u_{2},y-x\right\rangle
^{p}%
\]
and the points $x_{1}=y$, $x_{2}=x$, with the observation that
\[
\varphi\left(  -u_{1},-u_{2}\right)  =-\varphi\left(  u_{1},u_{2}\right)  .
\]
\hfill$\square$

\subsection{Localization\label{sectlocalizz}}

At the technical level, we do not need to localize the $\sigma$-finite
measures of the present work. However, we give a few remarks on localization
to help the intuitive interpretation of the model. Essentially we are going to
introduce rigorous analogs of the heuristic expressions (\ref{muformal}) and
(\ref{uformal}) written at the beginning of the previous section. The problem
there was that the law of $\xi^{\alpha}$ should be $\nu$, which is only a
$\sigma$-finite measure. For this reason one has to localize $\nu$ and $\mu$.

Given $A\in\mathcal{B}\left(  \Xi\right)  $ with $0<\nu\left(  A\right)
<\infty$, consider the measure $\mu_{A}$ defined as the restriction of $\mu$
to $A$:
\[
\mu_{A}\left(  B\right)  =\mu\left(  A\cap B\right)
\]
for any $B\in\mathcal{B}\left(  \Xi\right)  $. It can be written (the equality
is in law, or a.s. over a possibly enlarged probability space) as the sum of
independent random atoms each distributed according to the probability measure
$B\in\mathcal{B}\left(  A\right)  \mapsto\widetilde{\nu}_{A}\left(  B\right)
:=\nu(B|A)$:
\[
\mu_{A}(d\xi)=\sum_{\alpha=1}^{N_{A}}\delta_{\xi^{\alpha}}(d\xi)
\]
where $N_{A}$ is a Poisson random variable with intensity $\nu(A)$ and the
family of random variables $\{\xi^{\alpha}\}_{\alpha\in\mathbb{N}}$ is
independent and identically distributed according to $\widetilde{\nu}_{A}$.
Moreover if $\{A_{i}\}_{i\in\mathbb{N}}$ is a family of mutually disjoint
subsets of $\Xi$ then the r.v.s $\{\mu_{A_{i}}\}_{i\in\mathbb{N}}$ are independent.

Sets $A$ as above with a physical significance are the following ones. Given
$0<\eta<1$ and $R>0$, let
\[
A_{\eta,R}=\{(U,T,\ell,X)\in\Xi:\ell>\eta,|x_{0}|\leq R\}.
\]
In a fluid model we meet these sets if we consider only vortexes up to some
scale $\eta$ (it could be the Kolmogorov dissipation scale) and roughly
confined in a ball of radius $R$. If we assume that the measure $\gamma$
satisfies $0<\gamma(\ell>\eta)<\infty$ for each $\eta>0$, then $0<\nu\left(
A_{\eta,R}\right)  <\infty$, and for the measure $\mu_{\eta,R}:=\mu
_{A_{\eta,R}}$ we have the representation
\[
\mu_{\eta,R}(d\xi)=\sum_{\alpha=1}^{N_{\eta,R}}\delta_{\xi^{\alpha}}(d\xi)
\]
where $N_{\eta,R}$ is $\mathcal{P}\left(  \nu(A_{\eta,R})\right)  $ and
$\{\xi^{\alpha}\}_{\alpha\in\mathbb{N}}$ are i.i.d. with law $\widetilde{\nu
}_{\eta,R}:=\widetilde{\nu}_{A_{\eta,R}}$.

For any $x\in\mathbb{R}^{3}$ we may consider $\xi\mapsto u_{\text{single}%
}^{\xi}\left(  x\right)  $ as a random variable in $\mathbb{R}^{3}$, defined
$\widetilde{\mathcal{\nu}}_{\eta,R}$-a.s. on $A_{\eta,R}$. Moreover we may
consider the r.v. $u_{\eta,R}(x)$ on $\left(  \Omega,\mathcal{A},P\right)  $
defined as
\begin{equation}
u_{\eta,R}(x):=\int_{A_{\eta,R}}u_{\text{single}}^{\xi}\left(  x\right)
\mu(d\xi)=\int_{\Xi}u_{\text{single}}^{\xi}\left(  x\right)  \mu_{\eta,R}%
(d\xi).
\end{equation}
It is the velocity field at point $x$, generated by the vortex filaments in
$A_{\eta,R}\in\mathcal{B}\left(  \Xi\right)  $. In this case we have the
representation
\[
u_{\eta,R}(x)=\sum_{\alpha=1}^{N_{\eta,R}}u_{\text{single}}^{\xi^{\alpha}%
}\left(  x\right)  =\sum_{\alpha=1}^{N_{\eta,R}}\frac{U^{\alpha}}%
{(\ell^{\alpha})^{2}}\int_{0}^{T^{\alpha}}K_{\ell}(x-X_{t}^{\alpha})\wedge
dX_{t}^{\alpha}%
\]
where the quadruples $\xi^{\alpha}=(U^{\alpha},\ell^{\alpha},T^{\alpha
},X^{\alpha})$ are distributed according to $\widetilde{\nu}_{\eta,R}$ and are
independent. If $\omega\mapsto\xi\left(  \omega\right)  $ is any one of such
quadruples, the random variable
\[
\omega\mapsto u_{\text{single}}^{\xi\left(  \omega\right)  }\left(  x\right)
\]
is well defined, since the law of $\xi$ is $\widetilde{\nu}_{\eta,R}$ and the
random variable $\xi\mapsto u_{\text{single}}^{\xi}\left(  x\right)  $ is well
defined $\widetilde{\mathcal{\nu}}_{\eta,R}$-a.s. on $A_{\eta,R}$. Therefore
$u_{\eta,R}(x)$ is a well defined random variable on $\left(  \Omega
,\mathcal{A},P\right)  $. We have noticed this in contrast to the fact that
the definition of $u(x)$ required difficult estimates, because of the
contribution of infinitely many filaments.

Given the Poisson random field, by localization we have constructed the
velocity fields $u_{\eta,R}(x)$ that have a reasonable intuitive
interpretation. Connections between $u_{\eta,R}(x)$ and $u(x)$ can be
established rigorously at various levels. We limit ourselves to the following
example of statement.

\begin{lemma}
Assume
\[
\gamma(U^{p}\ell T)<\infty
\]
for any $p>0$. Then, at any $x\in\mathbb{R}^{3}$,
\[
\lim_{\left(  \eta,R^{-1}\right)  \rightarrow\left(  0^{+},0^{+}\right)
}E\left[  \left|  u_{\eta,R}(x)-u(x)\right|  ^{p}\right]  =0.
\]
\end{lemma}

\textsc{Proof.}\hspace{0.3cm} Since
\[
u_{\eta,R}(x)-u(x)=\int_{\Xi}\left(  1_{A_{\eta,R}}-1\right)  u_{\text{single}%
}^{\xi}\left(  x\right)  \mu(d\xi)
\]
we have
\begin{align*}
E\left[  \left|  u_{\eta,R}(x)-u(x)\right|  ^{p}\right]   &  \leq C_{p}%
\nu\left[  \left|  \left(  1_{A_{\eta,R}}-1\right)  u_{\text{single}}^{\cdot
}\left(  x\right)  \right|  ^{p}\right] \\
&  =C_{p}\nu\left[  \left(  1_{A_{\eta,R}}-1\right)  \left|  u_{\text{single}%
}^{\cdot}\left(  x\right)  \right|  ^{p}\right]
\end{align*}
Notice that we do not have $\nu\left(  A_{\eta,R}^{c}\right)  \rightarrow0$,
in general, so the argument to prove the lemma must take into account the
properties of the r.v. $u_{\text{single}}^{\cdot}\left(  x\right)  $. We have
(by Burkholder-Davis-Gundy inequality)
\begin{align*}
&  \nu\left[  \left(  1_{A_{\eta,R}}-1\right)  \left|  u_{\text{single}%
}^{\cdot}\left(  x\right)  \right|  ^{p}\right] \\
&  =\int_{\left|  x_{0}\right|  \geq R}dx_{0}\int_{l\geq\eta}\mathcal{W}%
_{x_{0}}(\left|  u_{\text{single}}^{\cdot}\left(  x\right)  \right|
^{p})d\gamma(U,\ell,T)\\
&  \leq C_{p}\int_{\left|  x_{0}\right|  \geq R}dx_{0}\int_{l\geq\eta}%
\frac{U^{p}}{\ell^{2p}}\mathcal{W}_{x_{0}}\left[  W(x)^{p/2}\right]
d\gamma(U,\ell,T)\\
&  \leq C_{p}\int\frac{U^{p}}{\ell^{2p}}d\gamma(U,\ell,T)\int_{\left|
x_{0}\right|  \geq R}\mathcal{W}_{x_{0}}\left[  W(x)^{p/2}\right]  dx_{0}%
\end{align*}
where $W(x)$ has been defined in lemma \ref{singlebound}. Let us show that
\begin{equation}
\int_{\left|  x_{0}\right|  \geq R}\mathcal{W}_{x_{0}}\left[  W(x)^{p/2}%
\right]  dx_{0}\leq C_{p}\ell^{2p}\ell T\cdot\theta\left(  R\right)
\label{key}%
\end{equation}
where $\theta\left(  R\right)  \rightarrow0$ as $R\rightarrow\infty$. The
proof will be complete after this result. Recall from the proof of lemma
\ref{singlebound} that we have
\begin{align*}
W(x)^{p/2}  &  \leq C_{p}\ell^{p}(L_{B(x,\ell)}^{T})^{p/2}+C_{p}\ell
^{p}\left(  \frac{\ell}{\Lambda}\right)  ^{2p}T^{p/2}\\
&  +C_{p}\ell^{p}\sum_{i=1}^{N-1}\left[  \left(  \frac{\ell}{\ell_{i-1}%
}\right)  ^{2}-\left(  \frac{\ell}{\ell_{i}}\right)  ^{2}\right]  \left(
\frac{\ell}{\ell_{i}}\right)  ^{p}(L_{B(x,\ell_{i})}^{T})^{p/2}.
\end{align*}
>From lemma \ref{lemmaultimo} we have
\begin{align*}
&  \int_{\left|  x_{0}\right|  \geq R}\mathcal{W}_{x_{0}}\left[  \left(
L_{B(x,\lambda)}^{T}\right)  ^{p/2}\right]  dx_{0}\\
&  \leq C_{p}(\lambda\wedge\sqrt{T})^{p}\lambda(\lambda\vee\sqrt{T})^{2}%
\exp\left(  -\frac{R-\left(  \left|  x\right|  +\lambda\right)  }{\sqrt{T}%
}\right)  .
\end{align*}
Hence
\begin{align*}
&  \int_{\left|  x_{0}\right|  \ \geq R}\mathcal{W}_{x_{0}}\left[
W(x)^{p/2}\right]  dx_{0}\\
&  \leq C_{p}\ell^{2p}\ell T\exp\left(  -\frac{R-\left(  \left|  x\right|
+\ell\right)  }{\sqrt{T}}\right)  +C_{p}\ell^{p}\left(  \frac{\ell}{\Lambda
}\right)  ^{2p}T^{p/2}\\
&  +C_{p}\ell^{p}\sum_{i=1}^{N-1}\left[  \left(  \frac{\ell}{\ell_{i-1}%
}\right)  ^{2}-\left(  \frac{\ell}{\ell_{i}}\right)  ^{2}\right]  \left(
\frac{\ell}{\ell_{i}}\right)  ^{p}\\
&  \times(\ell_{i}\wedge\sqrt{T})^{p}\ell_{i}(\ell_{i}\vee\sqrt{T})^{2}%
\exp\left(  -\frac{R-\left(  \left|  x\right|  +\ell_{i}\right)  }{\sqrt{T}%
}\right)  .
\end{align*}
Repeating the arguments of lemma \ref{singlebound} we arrive at
\begin{align*}
&  \int_{\left|  x_{0}\right|  \ \geq R}\mathcal{W}_{x_{0}}\left[
W(x)^{p/2}\right]  dx_{0}\\
&  \leq C_{p}\ell^{2p}\ell T\exp\left(  -\frac{R-\left(  \left|  x\right|
+\Lambda\right)  }{\sqrt{T}}\right)  +C_{p}\ell^{p}\left(  \frac{\ell}%
{\Lambda}\right)  ^{2p}T^{p/2}.
\end{align*}
Since $T$ and $\ell$ are smaller than one, and $p\geq2$, we also have
\begin{align*}
&  \int_{\left|  x_{0}\right|  \ \geq R}\mathcal{W}_{x_{0}}\left[
W(x)^{p/2}\right]  dx_{0}\\
&  \leq C_{p}\ell^{2p}\ell T\left(  \exp\left(  -R-\left(  \left|  x\right|
+\Lambda\right)  \right)  +\Lambda^{-2p}\right)  .
\end{align*}
With the choice $\Lambda=R/2$ we prove (\ref{key}). The proof is complete.
\hfill$\square$

\section{The structure function}

Given the random velocity field $u\left(  x,\cdot\right)  $ constructed above,
under the assumption of Corollary \ref{corollarybase}, a quantity of major
interest in the theory of turbulence is the \emph{longitudinal structure
function} defined for every integer $p$ and $\varepsilon>0$ as
\begin{equation}
S_{p}^{\parallel}(\varepsilon)=\expect[\langle e, u(x+\eps e) - u(x)\rangle
]^{p}%
\end{equation}
where $\langle\cdot,\cdot\rangle$ is the Euclidean scalar product in
$\mathbb{R}^{3}$, $e\in\mathbb{R}^{3}$ is a unit vector and $x\in
\mathbb{R}^{3}$, and $E$, we recall, is the expectation on $\left(
\Omega,\mathcal{A},P\right)  $.

\begin{remark}
\label{advert}We warn the reader familiar with the literature on statistical
fluid mechanics that $\varepsilon$ here is not the dissipation energy, but
just the spatial scale parameter. In the physical literature, it is commonly
denoted by $\ell$; however, in our mathematical analysis we need two
parameters: the scale parameter of the statistical observation, which we
denote by $\varepsilon$, and a parameter internal to the model that describes
the thickness of the different vortex filaments, that we denote by $\ell$.
\end{remark}

The moments $S_{p}^{\parallel}(\varepsilon)$ depend only on $\varepsilon$ and
$p$, since $u\left(  x,\cdot\right)  $ is homogeneous and isotropic: if
$e=R\cdot e_{1}$ where $e_{1}$ is a given unit vector, since the dual of $R$
is $R^{-1}$, we have
\begin{align*}
E\left[  \left\langle u(x+\varepsilon e)-u(x),e\right\rangle ^{p}\right]   &
=E\left[  \left\langle u(\varepsilon e)-u(0),e\right\rangle ^{p}\right] \\
&  =E\left[  \left\langle Ru(\varepsilon e_{1})-u(0),e\right\rangle
^{p}\right] \\
&  =E\left[  \left\langle u(\varepsilon e_{1})-u(0),e_{1}\right\rangle
^{p}\right]  .
\end{align*}
For this reason we do not write explicitly the dependence on $x$ and $e$.

Let us also introduce the \emph{(non-directional) structure function}
\[
S_{p}(\varepsilon)=E\left[  \left|  u(x+\varepsilon e)-u(x)\right|
^{p}\right]
\]
which also depends only on $\varepsilon$ and $p$. We obviously have
\[
S_{p}^{\parallel}(\varepsilon)\leq S_{p}(\varepsilon).
\]
We shall see that, for even integers $p$, they have the same scaling
properties. At the technical level, due to the previous inequality, it will be
sufficient to estimate carefully $S_{p}^{\parallel}(\varepsilon)$ from below
and $S_{p}(\varepsilon)$ from above.

\subsection{The main result}

The quantities $S_{p}^{\parallel}(\varepsilon)$ and $S_{p}(\varepsilon)$
describe the statistical behavior of the increments of the velocity field when
$\varepsilon\rightarrow0$ and have been extensively investigated, see
\cite{frisch}. They are expected to have a characteristic power-like behavior
of the form (\ref{scal}), and similarly for $S_{p}^{\parallel}(\varepsilon)$.
Our aim is to prove that, for the model described in the previous section with
a suitable choice of $\gamma$, (\ref{scal}) holds true in the sense that the
limit
\begin{equation}
\zeta_{p}=\lim_{\varepsilon\rightarrow0}\frac{\log S_{p}(\varepsilon)}%
{\log\varepsilon}%
\end{equation}
exists (similarly for $S_{p}^{\parallel}(\varepsilon)$) and is computable. The
following theorem gives us the necessary estimates from above and below, for a
rather general measure $\gamma$. Then, in the next subsection, we make a
choice of $\gamma$ in order to have the classical multifractal scaling theory.

\begin{theorem}
\label{maintheorem} Assume that
\[
\gamma(U^{p}\ell T)<\infty
\]
for every $p>1$. Then, for any even integer $p>1$ there exist two constants
$C_{p},c_{p}>0$ such that
\begin{equation}
c_{p}\gamma\left[  U^{p}\ell T1_{\ell<\varepsilon}\right]  \leq S_{p}%
^{\parallel}(\varepsilon)\leq S_{p}(\varepsilon)\leq C_{p}\gamma\left[
U^{p}\left(  \frac{\ell\wedge\varepsilon}{\ell}\right)  ^{p}\ell T\right]
\label{eq:general-scaling}%
\end{equation}
for every $\varepsilon\in\left(  0,1\right)  $.
\end{theorem}

The proof of this result is long and reported in a separate section.

We would like to give a very rough heuristic that could explain this result.
It must be said that we would not believe in this heuristic without the proof,
since some steps are too vague (we have devised this heuristic only a posteriori).

What we are going to explain is that
\[
\mathcal{W}(\left|  u_{\text{single}}\left(  x+\varepsilon e\right)
-u_{\text{single}}\left(  x\right)  \right|  ^{p})\sim U^{p}\left(  \frac
{\ell\wedge\varepsilon}{\ell}\right)  ^{p}\ell T.
\]
This is the hard part of the estimate.

Let us discuss separately the case $\varepsilon>\ell$ from the opposite one.
When $\varepsilon>\ell$ the vortex structure $u_{\text{single}}$ is very thin
compared to the length $\varepsilon$ of observation of the displacement, thus
the difference $u_{\text{single}}\left(  x+\varepsilon e\right)
-u_{\text{single}}\left(  x\right)  $ does not really play a role and the
value of $\mathcal{W}(\left|  u_{\text{single}}\left(  x+\varepsilon e\right)
-u_{\text{single}}\left(  x\right)  \right|  ^{p})$ comes roughly from the
separate contributions of $u_{\text{single}}\left(  x+\varepsilon e\right)  $
and $u_{\text{single}}\left(  x\right)  $, which are similar. Let us compute
$\mathcal{W}(\left|  u_{\text{single}}\left(  x\right)  \right|  ^{p})$.

Consider the expression (\ref{eq:single-vortex-u-ito}) which defines
$u_{\text{single}}(x)$. Strictly speaking, consider the short-range case,
otherwise there is a correction which makes even more difficult the intuition.
Very roughly, $K_{\ell}(x-X_{t})$ behaves like $\ell\cdot1_{X_{t}\in B\left(
x,\ell\right)  }$, hence, even more roughly, $u_{\text{single}}(x)$ behaves
like
\[
u_{\text{single}}(x)\sim\frac{U}{\ell}\int_{0}^{T}1_{X_{t}\in B\left(
x,\ell\right)  }dX_{t}.
\]
When $X_{t}$ is a smooth curve, say a stright line (at distances compared to
$\ell$), then $\int_{0}^{T}1_{X_{t}\in B\left(  x,\ell\right)  }dX_{t}$ is
roughly proportional to $\ell$ if $X_{t}$ crosses $B\left(  x,\ell\right)  $,
while it is zero otherwise. We assume the same result holds true when $X_{t}$
is a Brownian motion. In addition, $X_{t}$ crosses $B\left(  x,\ell\right)  $
with a probability proportional to the volume of the Wiener sausage, which is
$\ell T$. Summarizing, we have%
\[
\int_{0}^{T}1_{X_{t}\in B\left(  x,\ell\right)  }dX_{t}\sim\left\{
\begin{array}
[c]{cc}%
\ell & \text{with probability }\ell T\\
0 & \text{otherwise}%
\end{array}
\right.  .
\]
Therefore $u_{\text{single}}(x)$ takes rougly two values, $U$ with probability
$\ell T$ and $0$ otherwise. It follows that $\mathcal{W}(\left|
u_{\text{single}}\left(  x\right)  \right|  ^{p})\sim U^{p}\ell T$.

Consider now the case $\varepsilon<\ell$. The difference now is important.
Since the gradient of $K_{\ell}$ is of order one, we have
\[
u_{\text{single}}\left(  x+\varepsilon e\right)  -u_{\text{single}}\left(
x\right)  \sim\frac{U}{\ell^{2}}\varepsilon\int_{0}^{T}1_{X_{t}\in B\left(
x,\ell\right)  }dX_{t}.
\]
As above we conclude that $u_{\text{single}}\left(  x+\varepsilon e\right)
-u_{\text{single}}\left(  x\right)  $ takes rougly two values, $\varepsilon
{U}/{\ell}$ with probability $\ell T$ and $0$ otherwise. It follows that
$\mathcal{W}(\left|  u_{\text{single}}\left(  x\right)  \right|  ^{p})\sim
U^{p}\left(  \frac{\varepsilon}{\ell}\right)  ^{p}\ell T$. The intuitive
argument is complete.

\subsection{Example: the multifractal model\label{multifrexample}}

The most elementary idea to introduce a measure $\gamma$ on the parameters is
to take $U$ and $T$ as suitable powers of $\ell$, thus prescribing a relation
between the thickness $\ell$ and the length and intensity. That is, a relation
of the form
\[
d\gamma(U,\ell,T)=\delta_{\ell^{h}}(U)\delta_{\ell^{a}}(T)\ell^{-b}d\ell.
\]
Moreover, we have to prescribe the distribution of $\ell$ itself, which could
again be given by a power law $\ell^{-b}d\ell$. The K41 scaling described
below is such an example.

Having in mind multi-scale phenomena related to intermittency, we consider a
superposition of the previous scheme. Take a probability measure $\theta$ on
an interval $I\subset\mathbb{R}_{+}$ (which measures the relative importance
of the scaling exponent $h\in I$). Given two functions $a,b:I\rightarrow
\mathbb{R}_{+}$ with $a(h)\leq2$ (to ensure $\ell^{2}\leq T$) consider the
measure
\begin{equation}
d\gamma(U,\ell,T)=\int_{I}\left[  \delta_{\ell^{h}}(U)\delta_{\ell^{a(h)}%
}(T)\ell^{-b(h)}d\ell\right]  \theta(dh).
\end{equation}
Then, according to Theorem~\ref{maintheorem}
%
%
%
, we must evaluate
\begin{equation}%
\begin{split}
\gamma(U^{p}\ell T1_{\ell<\varepsilon})  &  =\int_{[0,\varepsilon]\times
I}\ell^{hp+1+a(h)-b(h)}d\ell\theta(dh)\\
&  =\int_{I}c_{p,h}\varepsilon^{hp+2+a(h)-b(h)}\theta(dh)
\end{split}
\label{estmultifr}%
\end{equation}
while
\[%
\begin{split}
\gamma\left[  U^{p}\left(  \frac{\varepsilon\wedge\ell}{\ell}\right)  ^{p}\ell
T\right]   &  =\int_{[0,\varepsilon]\times I}\ell^{hp+1+a(h)-b(h)}d\ell
\theta(dh)+\varepsilon^{p}\int_{[\varepsilon,+\infty)\times I}\ell
^{(h-1)p+1+a(h)-b(h)}d\ell\theta(dh)\\
&  =\int_{I}C_{p,h}\varepsilon^{hp+2+a(h)-b(h)}\theta(dh)
\end{split}
\]
As $\varepsilon\rightarrow0$, by Laplace method, we get
\[
\lim_{\varepsilon\rightarrow0}\frac{\log S_{p}(\varepsilon)}{\log\varepsilon
}=\inf_{h\in I}[hp+3-D(h)]=\zeta_{p}%
\]
with $D(h)=b(h)-a(h)+1$. With this choice of $\gamma$ we have recovered the
scaling properties of the \textit{multifractal model} of \cite{Parisi-Frisch}.
See \cite{frisch} for a review.

Consider the specific choice $\theta(dh)=\delta_{1/3}(dh)$ and $a(1/3)=2$,
$b(1/3)=4$. We have
\[
\zeta_{p}=[hp+3-D(h)]_{h=1/3}=\frac{p}{3}.
\]
This is the Kolmogorov K41 scaling law for 3d turbulence. The choice $a(1/3)=2
$, namely $T=\ell^{2}$, has the following geometrical meaning: the spatial
displacement and the thickness of the structure is comparable (remember that
the curves are Brownian), hence their shape is blob-like, as in the classical
discussions of ``eddies'' around K41. The choice $b(1/3)=4$, namely the
measure $\ell^{-4}d\ell$ for the parameter $\ell$, corresponds to the idea
that the eddies are space-filling: it is easy to see that in a box of unit
volume the number of eddies of size larger that $\ell$ is of order $\ell^{-3}%
$. Finally, the choice $\theta(dh)=\delta_{1/3}(dh)$, namely $U=\ell^{1/3}$,
is the key point that produces $\zeta_{p}=\frac{p}{3}$; one may attempt to
justify it by dimensional analysis or other means, but it essentially one of
the issues that should require a better foundation.

\begin{remark}
K41 can be obtain from this model also for other choices of the functions
$a(h)$ and $b(h)$. For example, we may take $\theta(dh)=\delta_{1/3}(dh)$,
$a(1/3)=0$ and $b(1/3)=2$. The value choice $\theta(dh)=\delta_{1/3}(dh)$ is
again the essential point. The value $a(1/3)=0$ means $T=1$, hence the model
is made of thin filaments of length comparable to the integral scale, instead
of blob-like objects. The value $b(1/3)=2$, namely the measure $\ell^{-2}%
d\ell$ for the parameter $\ell$, is again space filling in view of the major
length of the single vortex structures (recall also that the Hausdorff
dimension of Brownian trajectories is 2). From this example we see that in the
model described here it is possible to reconcile K41 with a geometry made of
thin long vortexes.
\end{remark}

\begin{remark}
\label{remmultmodif}Assume $\inf I>0$ and, for instance, $\sup_{I}D<\infty$.
Then $\lim_{p\rightarrow\infty}\zeta_{p}=+\infty$. In particular, $\zeta
_{p}>3$ for some even integer $p$. Since, by (\ref{estmultifr}) and Laplace
method,
\[
\gamma(U^{p}\ell T1_{\ell<\varepsilon^{1-\beta}})\leq C_{p,\beta}%
\varepsilon^{\zeta_{p}\left(  1-\beta\right)  }%
\]
for any $\beta>0$, the condition of remark \ref{remregularity} is satisfied.
Therefore the velocity field has a continuous modification.
\end{remark}

\section{Proof of Theorem \ref{maintheorem}}

Let us introduce some objects related to the structure functions at the level
of a single vortex filament. Let $e$ be a given unit vector, the first
canonical one to fix the ideas. Since $u(x)=\mu\left(  u_{\text{single}%
}^{\cdot}\left(  x\right)  \right)  $, then
\[
\left\langle u(\varepsilon e)-u(0),e\right\rangle =\mu\left[  \left\langle
\delta_{\varepsilon}u_{\text{single}},e\right\rangle \right]
\]
where
\[
\delta_{\varepsilon}u_{\text{single}}:=u_{\text{single}}^{\cdot}\left(
\varepsilon e\right)  -u_{\text{single}}^{\cdot}\left(  0\right)
\]
and similarly
\[
\left|  u(\varepsilon e)-u(0)\right|  =\left|  \mu\left[  \delta_{\varepsilon
}u_{\text{single}}\right]  \right|  .
\]
Of major technical interest will be the quantities, of structure function
type,
\[
\mathcal{S}_{p}^{e}(\varepsilon)=\mathcal{W}(\left\langle \delta_{\varepsilon
}u_{\text{single}},e\right\rangle ^{p})
\]%
\[
\mathcal{S}_{p}(\varepsilon)=\mathcal{W}(\left|  \delta_{\varepsilon
}u_{\text{single}}\right|  ^{p}).
\]
They depend also on $\ell,T,U$.

\subsection{Lower bound}

As a direct consequence of lemma \ref{reflectsymm} and lemma \ref{moments} we have:

\begin{corollary}
If $k$ is an odd positive integer, then
\[
S_{k}^{\parallel}(\varepsilon)=0.
\]
Moreover, for any even integer $p>1$ there is a constant $c_{p}$ such that
\begin{align*}
S_{p}^{\parallel}(\varepsilon)  &  \geq c_{p}\nu\left[  \left\langle
\delta_{\varepsilon}u_{\text{single}},e\right\rangle ^{p}\right] \\
&  =c_{p}\gamma\left[  \mathcal{S}_{p}^{e}(\varepsilon)\right]  .
\end{align*}
\end{corollary}

If we prove that
\[
\mathcal{S}_{p}^{e}(\varepsilon)\geq c_{p}^{\prime}U^{p}\ell T
\]
for every $\varepsilon\in\left(  \ell,1\right)  $ and some constant
$c_{p}^{\prime}>0$, then
\[
\gamma\left[  \mathcal{S}_{p}^{e}(\varepsilon)\right]  \geq c_{p}%
^{\prime\prime}\gamma\left[  U^{p}\ell T1_{\ell<\varepsilon}\right]
\]
for every $\varepsilon\in\left(  0,1\right)  $ and some constant
$c_{p}^{\prime\prime}>0$, which implies the lower bound of theorem
\ref{maintheorem}. We have
\[
\mathcal{S}_{p}^{e}(\varepsilon)=\int\mathcal{W}_{x_{0}}\left[  \left\langle
\delta_{\varepsilon}u_{\text{single}},e\right\rangle ^{p}\right]  dx_{0}.
\]
Since
\[
\left\langle \delta_{\varepsilon}u_{\text{single}},e\right\rangle =\frac
{U}{\ell^{2}}\int_{0}^{T}\left\langle K_{\ell}^{e}(\varepsilon e-X_{t}%
)-K_{\ell}^{e}(0-X_{t}),dX_{t}\right\rangle
\]
where%
\[
K_{\ell}^{e}\left(  y\right)  =K_{\ell}\left(  y\right)  \wedge e
\]
by Burkholder-Davis-Gundy inequality, there is $c_{p}>0$ such that
\[
\mathcal{W}_{x_{0}}\left[  \left\langle \delta_{\varepsilon}u_{\text{single}%
},e\right\rangle ^{p}\right]  \geq c_{p}\frac{U^{p}}{\ell^{2p}}\mathcal{W}%
_{x_{0}}\left[  \left(  W_{\varepsilon}^{e}\right)  ^{p/2}\right]
\]
where
\[
W_{\varepsilon}^{e}=\int_{0}^{T}dt|K_{\ell}^{e}(\varepsilon e-X_{t})-K_{\ell
}^{e}(0-X_{t})|^{2}.
\]
Here $K_{\ell}^{e}(y)=\left\langle K_{\ell}(y),e\right\rangle $. Therefore
\begin{align*}
\mathcal{S}_{p}^{e}(\varepsilon)  &  \geq c_{p}\frac{U^{p}}{\ell^{2p}}%
\int\mathcal{W}_{x_{0}}\left[  \left(  W_{\varepsilon}^{e}\right)
^{p/2}\right]  dx_{0}\\
&  =c_{p}\frac{U^{p}}{\ell^{2p}}\mathcal{W}\left[  \left(  W_{\varepsilon}%
^{e}\right)  ^{p/2}\right]  .
\end{align*}
The proof of the theorem is then complete with the following basic estimate.

\begin{lemma}
Given $p>0$, there exist $c_{p}>0$ such that for every $\ell^{2}\leq T\leq1$
and
\[
\varepsilon>\ell
\]
we have
\[
\mathcal{W}\left[  \left(  W_{\varepsilon}^{e}\right)  ^{p/2}\right]  \geq
c_{p}\ell^{2p}\ell T.
\]
\end{lemma}

\textsc{Proof.}\hspace{0.3cm}Recall that
\[
K_{1}(x)=-2Q\frac{x}{|x|^{3}}\quad\text{for $|x|\geq$}1.
\]
Consider the function
\begin{align*}
f\left(  z,\alpha\right)   &  =K_{1}^{e}(z)-K_{1}^{e}(z+\alpha e)\\
&  =\left(  K_{1}(z)-K_{1}(z+\alpha e)\right)  \wedge e
\end{align*}
defined for $z\in\mathbb{R}^{3}$ and $\alpha\geq1$. Let $e^{\perp}$ be any
unit vector orthogonal to $e$. We have
\begin{align*}
\left|  f\left(  e^{\perp},\alpha\right)  \right|   &  =2Q\left|
1-\frac{\left(  e^{\perp}+\alpha e\right)  \wedge e}{\left|  e^{\perp}+\alpha
e\right|  ^{3}}\right| \\
&  =2Q\left|  1-\frac{1}{\left|  e^{\perp}+\alpha e\right|  ^{3}}\right| \\
&  \geq f\left(  e^{\perp},1\right)  =C_{0}Q
\end{align*}
for $C_{0}=2\left(  1-2^{-\frac{3}{2}}\right)  $. Moreover, $K_{1}$ is
globally Lipschitz, hence there is $L>0$ such that
\begin{align*}
&  \left|  f\left(  e^{\perp}+w,\alpha\right)  -f\left(  e^{\perp}%
,\alpha\right)  \right| \\
&  \leq\left|  K_{1}(e^{\perp}+w)-K_{1}(e^{\perp})\right|  +\left|
K_{1}(e^{\perp}+\alpha e+w)-K_{1}(e^{\perp}+\alpha e)\right|  \leq L\left|
w\right|
\end{align*}
for every $w\in\mathbb{R}^{3}$ and $\alpha\geq1$. Therefore, there exists a
ball $B(e^{\perp},a)\subset B(0,2)$ and a constant $C_{1}>0$ such that when
$z\in B(e^{\perp},a)$ we have
\[
|K_{1}^{e}(z)-K_{1}^{e}(z+\alpha e)|\geq C_{1}%
\]
uniformly in $\alpha\geq1$. Then reintroducing the scaling factor $\ell$ we
obtain that for $y\in B(\ell e^{\perp},\ell a)$
\[
|K_{\ell}^{e}(y)-K_{\ell}^{e}(y+\varepsilon e)|>C_{1}\ell
\]
uniformly in $\ell\in(0,1)$ and $\varepsilon>\ell$. Then we have
\[
|K_{\ell}^{e}(0-X_{t})-K_{\ell}^{e}(\varepsilon e-X_{t})|\geq C_{1}%
\ell1_{X_{t}\in B(-\ell e^{\perp},\ell a)}.
\]
Hence, if $\varepsilon>\ell$ we have
\begin{align*}
W_{\varepsilon}^{e}  &  =\int_{0}^{T}dt|K_{\ell}^{e}(0-X_{t})-K_{\ell}%
^{e}(\varepsilon e-X_{t})|^{2}\\
&  \geq C_{1}\ell^{2}L_{B(-\ell e^{\perp},\ell a)}^{T}.
\end{align*}

>From the lower bound in (\ref{eq:Lmoments-bound4}) proved in the next
section,
\[
\mathcal{W}\left[  \left(  W_{\varepsilon}^{e}\right)  ^{p/2}\right]  \geq
c_{p}\ell^{p}\mathcal{W}\left[  \left(  L_{B(-\ell e^{\perp},\ell a)}%
^{T}\right)  ^{p/2}\right]  \geq c_{p}^{\prime}\ell^{2p}\ell T.
\]
The lower bound is proved.

\begin{remark}
With a bit more of effort is is also possible to prove the bound
\[
\mathcal{S}_{p}^{e}(\varepsilon)\geq c_{p}U^{p}\left(  \frac{\varepsilon
\wedge\ell}{\ell}\right)  ^{p}\ell T
\]
valid for every $\varepsilon\in\left(  0,1\right)  $ (not only for
$\varepsilon>\ell$). This would be the same as the upper bound, but we do not
need it to prove that the behaviors as $\varepsilon\rightarrow0$ of the upper
and lower bound is the same.
\end{remark}

\subsection{Upper bound}

\begin{lemma}
For every even $p$ there exists a constant $C_{p}>0$ such that
\begin{align*}
S_{p}(\varepsilon)\leq C_{p}\nu\left[  \left|  \delta_{\varepsilon
}u_{\text{single}}\right|  ^{p}\right] \\
=C_{p}\gamma\left[  \mathcal{S}_{p}(\varepsilon)\right]  .
\end{align*}
\end{lemma}

\textsc{Proof.}\hspace{0.3cm} Let $\left[  \delta_{\varepsilon}%
u_{\text{single}}\left(  \xi\right)  \right]  _{i}$ be the $i$-th component of
$\delta_{\varepsilon}u_{\text{single}}\left(  \xi\right)  $. We have
\[
S_{p}(\varepsilon)\leq C_{p}\sum_{i=1}^{3}E\left[  \left(  \mu\left[  \left[
\delta_{\varepsilon}u_{\text{single}}\left(  \xi\right)  \right]  _{i}\right]
\right)  ^{p}\right]
\]
and thus, by lemma \ref{moments},
\[
S_{p}(\varepsilon)\leq C_{p}^{\prime}\sum_{i=1}^{3}\nu\left[  \left[
\delta_{\varepsilon}u_{\text{single}}\left(  \xi\right)  \right]  _{i}%
^{p}\right]
\]
which implies the claim. \hfill$\square$

It is then sufficient to prove the bound
\[
\mathcal{S}_{p}(\varepsilon)\leq C_{p}U^{p}\left(  \frac{\varepsilon\wedge
\ell}{\ell}\right)  ^{p}\ell T.
\]
Again as above, We have
\[
\mathcal{S}_{p}(\varepsilon)=\int\mathcal{W}_{x_{0}}\left[  \left|
\delta_{\varepsilon}u_{\text{single}}\right|  ^{p}\right]  dx_{0}%
\]
where, by Burkholder-Davis-Gundy inequality, there is $C_{p}>0$ such that
\[
\mathcal{W}_{x_{0}}\left[  \left|  \delta_{\varepsilon}u_{\text{single}%
}\right|  ^{p}\right]  \leq C_{p}\frac{U^{p}}{\ell^{2p}}\mathcal{W}_{x_{0}%
}\left[  W_{\varepsilon}^{p/2}\right]
\]
where
\[
W_{\varepsilon}=\int_{0}^{T}dt|K_{\ell}(\varepsilon e-X_{t})-K_{\ell}%
(0-X_{t})|^{2}.
\]
Therefore
\begin{align*}
\mathcal{S}_{p}(\varepsilon)  &  \leq C_{p}\frac{U^{p}}{\ell^{2p}}%
\int\mathcal{W}_{x_{0}}\left[  W_{\varepsilon}^{p/2}\right]  dx_{0}\\
&  =C_{p}\frac{U^{p}}{\ell^{2p}}\mathcal{W}\left[  W_{\varepsilon}%
^{p/2}\right]  .
\end{align*}
It is then sufficient to prove the bound
\[
\mathcal{W}\left[  W_{\varepsilon}^{p/2}\right]  \leq C_{p}\ell^{2p}\left(
\frac{\varepsilon\wedge\ell}{\ell}\right)  ^{p}\ell T.
\]

For $\varepsilon>\ell$ it is not necessary to keep into account the closeness
of $\varepsilon e$ to $0$: each term in the difference of $W_{\varepsilon}$
has already the necessary scaling. The hard part of the work has been done in
lemma \ref{singlebound} above.

\begin{lemma}
Given $p>0$, there exist $C_{p}>0$ such that for every $\ell^{2}\leq T\leq1$
and
\[
\varepsilon>\ell
\]
we have
\[
\mathcal{W}\left[  W_{\varepsilon}^{p/2}\right]  \leq C_{p}\ell^{2p}\ell T.
\]
\end{lemma}

\textsc{Proof.}\hspace{0.3cm} Since
\[
W_{\varepsilon}\leq2W(\varepsilon e)+2W(0)
\]
where $W(x)=\int_{0}^{T}\left|  K_{\ell}(x-X_{t})\right|  ^{2}dt$, by lemma
\ref{singlebound} we have the result. \hfill$\square$

For $\varepsilon\leq\ell$ need to extract a power of $\varepsilon$ from the
estimate of $\mathcal{W}\left[  W_{\varepsilon}^{p/2}\right]  $. We
essentially repeat the multi-scale argument in the proof of lemma
\ref{singlebound}, with suitable modifications.

\begin{lemma}
As in the previous lemma, when
\[
\varepsilon\leq\ell
\]
we have
\[
\mathcal{W}\left[  W_{\varepsilon}^{p/2}\right]  \leq C_{p}\varepsilon^{p}%
\ell^{p}\ell T.
\]
\end{lemma}

\textsc{Proof.}\hspace{0.3cm} Since now $\varepsilon$ is smaller that $\ell$
we bound $K_{\ell}(y)-K_{\ell}(z)$ for $|y-z|\leq\varepsilon$ as
\[
|K_{\ell}(y)-K_{\ell}(z)|\leq C\varepsilon
\]
if $|y|\leq2\ell$. If $|y|>2\ell$ then $|z|\geq\ell$ and using the explicit
form of the kernel $K_{\ell}$ we have the bound
\[
|K_{\ell}(y)-K_{\ell}(z)|\leq C\varepsilon\left(  \frac{\ell}{u}\right)
^{3}\leq C\varepsilon\left(  \frac{\ell}{u}\right)  ^{2}%
\]
where $u$ is the minimum between $|y|$ and $|z|$ and in this case $\ell/u<1$.
Then, given a partition of $\left[  2\ell,\Lambda\right]  $, say $2\ell
=\ell_{0}<\ell_{1}<...<\ell_{N}=\Lambda$, as in the proof of lemma
\ref{singlebound}, we get
\begin{align*}
\left|  K_{\ell}(y)-K_{\ell}(z)\right|  ^{2}  &  =\left|  K_{\ell}(y)-K_{\ell
}(z)\right|  ^{2}\left(  1_{|y|\leq2\ell}+\sum_{i=1}^{N}1_{\ell_{i-1}%
<|y|<\ell_{i}}+1_{|y|>\Lambda}\right) \\
&  \leq C\varepsilon^{2}1_{|y|\leq2\ell}+C\sum_{i=1}^{N}\varepsilon^{2}\left(
\frac{\ell}{\ell_{i-1}-\ell}\right)  ^{4}1_{\ell_{i-1}<|y|<\ell_{i}%
}+C\varepsilon^{2}\left(  \frac{\ell}{\Lambda}\right)  ^{4}1_{|y|>\Lambda}%
\end{align*}%
\[
\leq C\varepsilon^{2}1_{|y|\leq2\ell}+2^{4}C\sum_{i=1}^{N}\varepsilon
^{2}\left(  \frac{\ell}{\ell_{i-1}}\right)  ^{4}1_{\ell_{i-1}<|y|<\ell_{i}%
}+C\varepsilon^{2}\left(  \frac{\ell}{\Lambda}\right)  ^{4}1_{|y|>\Lambda}%
\]
where we have used the fact that $(u-\ell)^{-1}\leq{2}/{u}$ for $u\geq2\ell$.
Then
\[
W_{\varepsilon}\leq C\varepsilon^{2}L_{B(0,2\ell)}^{T}+C\sum_{i=1}%
^{N}\varepsilon^{2}\left(  \frac{\ell}{\ell_{i-1}}\right)  ^{4}L_{B(0,\ell
_{i})\backslash B(x,\ell_{i-1})}^{T}+C\varepsilon^{2}\left(  \frac{\ell
}{\Lambda}\right)  ^{4}L_{B(x,\Lambda)^{c}}^{T}%
\]
and arguing as in that proof
\begin{align*}
W_{\varepsilon}{}^{p/2}  &  \leq C\varepsilon^{p}(L_{B(0,2\ell)}^{T}%
)^{p/2}+C\varepsilon^{p}\left(  \frac{\ell}{\Lambda}\right)  ^{2p}%
[L_{B(0,\Lambda)}^{T}]^{p/2}+C\varepsilon^{p}\left(  \frac{\ell}{\Lambda
}\right)  ^{2p}\left(  L_{B(x,\Lambda)^{c}}^{T}\right)  ^{p/2}\\
&  +C\varepsilon^{p}\sum_{i=1}^{N-1}\left[  \left(  \frac{\ell}{\ell_{i-1}%
}\right)  ^{2}-\left(  \frac{\ell}{\ell_{i}}\right)  ^{2}\right]  \left(
\frac{\ell}{\ell_{i}}\right)  ^{p}(L_{B(0,\ell_{i})}^{T})^{p/2}.
\end{align*}
Then, from lemma \ref{lemma:definitive-L-bounds} in the form
\[
\mathcal{W}\left[  \left(  L_{B(x,\lambda)}^{T}\right)  ^{p/2}\right]  \leq
C_{p}(\lambda\wedge\sqrt{T})^{p}\lambda(\lambda\vee\sqrt{T})^{2}.
\]
and the obvious bound $L_{B}^{T}\leq T$ for every Borel set $B$, we have
\[
\mathcal{W}\left[  W_{\varepsilon}^{p/2}\right]  \leq C\varepsilon^{p}\ell
^{p}\ell T+C\varepsilon^{p}\left(  \frac{\ell}{\Lambda}\right)  ^{2p}T^{p/2}%
\]%

\[
+C\varepsilon^{p}\sum_{i=1}^{N-1}\left[  \left(  \frac{\ell}{\ell_{i-1}%
}\right)  ^{2}-\left(  \frac{\ell}{\ell_{i}}\right)  ^{2}\right]  \left(
\frac{\ell}{\ell_{i}}\right)  ^{p}(\ell_{i}\wedge\sqrt{T})^{p}\ell_{i}%
(\ell_{i}\vee\sqrt{T})^{2}%
\]
and taking the limit as the partition gets finer:
\[
\mathcal{W}\left[  W_{\varepsilon}^{p/2}\right]  \leq C\varepsilon^{p}\ell
^{p}\ell T+C\varepsilon^{p}\left(  \frac{\ell}{\Lambda}\right)  ^{2p}T^{p/2}%
\]%
\[
+C\varepsilon^{p}\int_{2\ell}^{\Lambda}\left(  \frac{\ell}{u}\right)
^{p}(u\wedge\sqrt{T})^{p}(u\vee\sqrt{T})^{2}ud\left[  -\left(  \frac{\ell}%
{u}\right)  ^{2}\right]  .
\]
A direct computation of the integral as in the proof of lemma
\ref{singlebound} completes the proof, taking the limit as $\Lambda
\rightarrow\infty$.\hfill$\square$

\section{Auxiliary results on Brownian occupation measure}

In this section we prove the estimates on $\mathcal{W}[|L_{B(u,\ell)}%
^{T}|^{p/2}]$ which constitute the technical core of the previous sections.
The literature on Brownian occupation measure is wide, so it is possible that
results proved here are given somewhere or may be deduced from known results.
However, we have not found the uniform estimates we need, so we prefer to give
a full self-contained proofs for completeness. Of course, several ideas we use
are inspired by the existing literature (in particular, a main source of
inspiration has been \cite{LGall}).

First, notice that $\mathcal{W}[|L_{B(u,\ell)}^{T}|^{p/2}]$ does not depend on
$u$. But this is not of great help. When $p=2$, $\mathcal{W}[|L_{B(u,\ell
)}^{T}|^{p/2}]$ can be explicitly computed:%

\begin{align*}
&  \mathcal{W}[L_{B(u,\ell)}^{T}]=\int_{\mathbb{R}^{3}}dx_{0}\int_{0}%
^{T}dt\int_{B(u,\ell)}dz\;p_{t}\left(  z-x_{0}\right) \\
&  =\int_{\mathbb{R}^{3}}dx_{0}\int_{0}^{T}dt\int_{B(u,\ell)}dz\;p_{t}\left(
x_{0}\right) \\
&  =|B(u,\ell)|\int_{0}^{T}dt\int_{\mathbb{R}^{3}}p_{t}\left(  x_{0}\right)
dx_{0}=\varpi\ell^{3}T
\end{align*}
where $\varpi$ is a geometrical constant and $p_{t}\left(  x\right)  $ is the
density of the 3D Brownian motion at time $t$. The estimate of $\mathcal{W}%
[|L_{B(u,\ell)}^{T}|^{p/2}]$ for general $p$ requires much more work.

Let
\[
\tau_{B(u,\ell)}=\inf\{t\geq0:X_{t}\in B(u,\ell)\}
\]
the entrance time in $B(u,\ell)$ for the canonical process. We continue to
denote by $\mathcal{W}_{x_{0}}$ the Wiener measure starting at $x_{0}$ and
also the mean value with respect to it;\ similarly for $\mathcal{W}$, the
$\sigma$-finite measure $d\mathcal{W}_{x_{0}}(X)dx_{0}$.

\begin{lemma}
For any $p>0$, $T>0$, $\ell>0$, $x_{0},u\in\mathbb{R}^{3}$, it holds that
\[
\mathcal{W}_{x_{0}}[\tau_{B(u,\ell/2)}\leq T/2]\mathcal{W}_{0}[|L_{B(0,\ell
/2)}^{T/2}|^{p}]\leq\mathcal{W}_{x_{0}}[|L_{B(u,\ell)}^{T}|^{p}]
\]%
\begin{equation}
\leq\mathcal{W}_{x_{0}}[\tau_{B(u,\ell)}\leq T]\mathcal{W}_{0}[|L_{B(0,2\ell
)}^{T}|^{p}]. \label{eq:Lmoments-bound}%
\end{equation}
\end{lemma}

\textsc{Proof.}\hspace{0.3cm} Let us prove the upper bound. Let us set, for
simplicity, $\tau=\tau_{B(u,\ell)}\wedge T$. When $\tau\leq t<T$ we have
\[
X_{t}\in B(u,\ell)\Rightarrow X_{t}-X_{\tau}\in B(0,2\ell)
\]
then
\[
L_{B(u,\ell)}^{T}=\int_{\tau}^{T}1_{X_{t}\in B(u,\ell)}dt\leq\int_{\tau}%
^{T}1_{\{(X_{t}-X_{\tau})\in B(0,2\ell)\}}dt=\int_{0}^{T-\tau}1_{(X_{\tau
+t}-X_{\tau})\in B(0,2\ell)}dt
\]
Then, taking into account that $\tau=T$ implies $L_{B(u,\ell)}^{T}=0$,
\[
L_{B(u,\ell)}^{T}\leq1_{\{\tau<T\}}\int_{0}^{T}1_{(X_{\tau+t}-X_{\tau})\in
B(0,2\ell)}dt
\]
which gives us, using the strong Markov property,
\[%
\begin{split}
\mathcal{W}_{x_{0}}|L_{B(u,\ell)}^{T}|^{p}  &  \leq\mathcal{W}_{x_{0}}\left[
1_{\tau<T}\left(  \int_{0}^{T}1_{(X_{\tau+t}-X_{\tau})\in B(0,2\ell
)}dt\right)  ^{p}\right] \\
&  =\mathcal{W}_{x_{0}}[\tau<T]\mathcal{W}_{0}\left[  \left(  \int_{0}%
^{T}1_{X_{t}\in B(0,2\ell)}dt\right)  ^{p}\right] \\
&  =\mathcal{W}_{x_{0}}[\tau_{B(u,\ell)}<T]\mathcal{W}_{0}\left[
(L_{B(0,2\ell)}^{T})^{p}\right]  .
\end{split}
\]
The upper bound is proved.

Let us proceed with the lower bound. Let $\tau^{\prime}=\tau_{B(u,\ell
/2)}\wedge T$. When $\tau^{\prime}\leq t\leq T$ we have
\[
X_{t}-X_{\tau^{\prime}}\in B(0,\ell/2)\Rightarrow X_{t}\in B(u,\ell)
\]
then
\[%
\begin{split}
L_{B(u,\ell)}^{T}  &  \geq\int_{\tau^{\prime}}^{T}1_{X_{t}\in B(u,\ell)}%
dt\geq\int_{\tau^{\prime}}^{T}1_{X_{t}-X_{\tau^{\prime}}\in B(0,\ell/2)}dt\\
&  \geq\int_{0}^{T-\tau^{\prime}}1_{X_{\tau^{\prime}+t}-X_{\tau^{\prime}}\in
B(0,\ell/2)}dt\\
&  \geq1_{\tau^{\prime}\leq T/2}\int_{0}^{T-\tau^{\prime}}1_{X_{\tau^{\prime
}+t}-X_{\tau^{\prime}}\in B(0,\ell/2)}dt\\
&  \geq1_{\tau^{\prime}\leq T/2}\int_{0}^{T/2}1_{X_{\tau^{\prime}+t}%
-X_{\tau^{\prime}}\in B(0,\ell/2)}dt.
\end{split}
\]
Then, using again the strong Markov property with respect to the stopping time
$\tau^{\prime}$, we obtain
\[
\mathcal{W}_{x_{0}}[|L_{B(u,\ell)}^{T}|^{p}]\geq\mathcal{W}_{x_{0}}%
[\tau_{B(u,\ell/2)}\leq T/2]\mathcal{W}_{0}[(L_{B(0,\ell/2)}^{T/2})^{p}]
\]
and the proof is complete. \hfill$\square$\\[0.5cm]

Letting $p=1$ in the previous lemma and using the scale invariance of BM we
obtain an upper bounds for $\mathcal{W}_{x_{0}}[\tau_{B(u,\ell)}\leq T]$ as
\[
\mathcal{W}_{x_{0}}[\tau_{B(u,\ell)}\leq T]=\mathcal{W}_{x_{0}/\sqrt{2}}%
[\tau_{B(u/\sqrt{2},\ell/\sqrt{2})}\leq T/2]\leq\frac{\mathcal{W}_{x_{0}%
/\sqrt{2}}[L_{B(u/\sqrt{2},\sqrt{2}\ell)}^{T}]}{\mathcal{W}_{0}[L_{B(0,\ell
/\sqrt{2})}^{T/2}]}%
\]
and the corresponding lower bound for $\mathcal{W}_{x_{0}}[\tau_{B(u,\ell
/2)}\leq T/2]$:
\[
\mathcal{W}_{x_{0}}[\tau_{B(u,\ell/2)}\leq T/2]=\mathcal{W}_{\sqrt{2}x_{0}%
}[\tau_{B(\sqrt{2}u,\ell/\sqrt{2})}\leq T]\geq\frac{\mathcal{W}_{\sqrt{2}%
x_{0}}[L_{B(\sqrt{2}u,\ell/\sqrt{2})}^{T}]}{\mathcal{W}_{0}[L_{B(0,\sqrt
{2}\ell)}^{T}]}%
\]
which leads to the easy corollary:

\begin{corollary}
For any $p>0$, $T>0$, $\ell>0$, $x_{0},u\in\mathbb{R}^{3}$, it holds that
\begin{multline}
\frac{\mathcal{W}_{\sqrt{2}x_{0}}[L_{B(\sqrt{2}u,\ell/\sqrt{2})}%
^{T}]\mathcal{W}_{0}[|L_{B(0,\ell/2)}^{T/2}|^{p}]}{\mathcal{W}_{0}%
[L_{B(0,\sqrt{2}\ell)}^{T/2}]}\leq\mathcal{W}_{x_{0}}[|L_{B(u,\ell)}^{T}%
|^{p}]\label{eq:Lmoments-bound2}\\
\leq\frac{\mathcal{W}_{x_{0}/\sqrt{2}}[L_{B(u/\sqrt{2},\sqrt{2}\ell)}%
^{T}]\mathcal{W}_{0}[|L_{B(0,2\ell)}^{T}|^{p}]}{\mathcal{W}_{0}[L_{B(0,\ell
/\sqrt{2})}^{T}]}%
\end{multline}
\end{corollary}

\begin{lemma}
\label{lemma:Lzero-moments} Given $\alpha>0$ and $p>0$, there exist constants
$c,C>0$ such that the following properties hold true: for every $T,\ell$
satisfying $T/\ell^{2}\geq\alpha$ we have
\[
c\ell^{p}\leq\mathcal{W}_{0}|L_{B(0,\ell)}^{T}|^{p/2}\leq C\ell^{p}%
\]
while if $T/\ell^{2}\leq\alpha$
\[
cT^{p/2}\leq\mathcal{W}_{0}|L_{B(0,\ell)}^{T}|^{p/2}\leq CT^{p/2}.
\]
\end{lemma}

\textsc{Proof.}\hspace{0.3cm} Consider first $T/\ell^{2}\geq\alpha$. In
distribution
\begin{equation}
L_{B(0,\ell)}^{T}\overset{\mathcal{L}}{=}\ell^{2}\int_{0}^{T/\ell^{2}}%
1_{X_{t}\in B(0,1)}dt \label{star}%
\end{equation}
and moreover we have that
\begin{equation}
c\leq\mathcal{W}_{0}\left[  \int_{0}^{T/\ell^{2}}1_{X_{t}\in B(0,1)}dt\right]
^{p}\leq C \label{starstar}%
\end{equation}
uniformly in $T/\ell^{2}\geq\alpha>0$ (the constants depend on $\alpha$ and
$p$). The lower bound is obtained by setting $T/\ell^{2}=\alpha$ while the
upper bound is given by the fact that
\[
\mathcal{W}_{0}\left[  \int_{0}^{\infty}1_{X_{t}\in B(0,1)}dt\right]
^{p}<\infty
\]
for any $p>0$ (see for instance \cite{Kac}, section 3). From (\ref{star}) and
(\ref{starstar}) we get the first claim of the lemma.

Next, if $T/\ell^{2}\leq\alpha$, in distribution:
\[
L_{B(0,\ell)}^{T}\overset{\mathcal{L}}{=}T\int_{0}^{1}1_{X_{t}\in
B(0,\ell/\sqrt{T})}dt
\]
and
\[
\mathcal{W}_{0}|L_{B(0,1/\sqrt{\alpha})}^{1}|^{p}\leq\mathcal{W}%
_{0}|L_{B(0,\ell/\sqrt{T})}^{1}|^{p}\leq1
\]
so the second claim is also proved.\hfill$\square$

\begin{lemma}
\label{lemma:definitive-L-bounds} Given $p>0$, $\alpha>0$, there are constants
$c,C>0$ such that, for $T\geq\alpha\ell^{2}$
\begin{equation}
c\ell^{p}\ell T\leq\mathcal{W}[|L_{B(u,\ell)}^{T}|^{p/2}]\leq C\ell^{p}\ell T
\label{eq:Lmoments-bound4}%
\end{equation}
while for $T\leq\alpha\varepsilon^{2}$
\begin{equation}
cT^{p/2}\varepsilon^{3}\leq\mathcal{W}[|L_{B(u,\varepsilon)}^{T}|^{p/2}]\leq
CT^{p/2}\varepsilon^{3}. \label{eq:Lmoments-bound4-bis}%
\end{equation}
\end{lemma}

\textsc{Proof.}\hspace{0.3cm}Let us prove (\ref{eq:Lmoments-bound4}). Using
lemma~\ref{lemma:Lzero-moments}, equation (\ref{eq:Lmoments-bound2}) becomes
\begin{equation}
c\ell^{p-2}\mathcal{W}_{\sqrt{2}x_{0}}[L_{B(\sqrt{2}u,\ell/\sqrt{2})}^{T}%
]\leq\mathcal{W}_{x_{0}}[|L_{B(u,\ell)}^{T}|^{p/2}]\leq C\ell^{p-2}%
\mathcal{W}_{x_{0}/\sqrt{2}}[L_{B(u/\sqrt{2},\sqrt{2}\ell)}^{T}] \label{dis1}%
\end{equation}
for two constants $c,C>0$ (depending on $\alpha$ and $p$). Moreover, as we
remarked at the beginning of the section,
\[
\mathcal{W}[L_{B(u,\ell)}^{T}]=\varpi\ell^{3}T.
\]
Using this identity in (\ref{dis1}), we get (\ref{eq:Lmoments-bound4}).

Now consider eq.(\ref{eq:Lmoments-bound4-bis}). Assume $T\leq\alpha
\varepsilon^{2}$. Using eq.(\ref{eq:Lmoments-bound2}) and
lemma~\ref{lemma:Lzero-moments} we have
\begin{equation}
cT^{\frac{p}{2}-1}\mathcal{W}_{\sqrt{2}x_{0}}[L_{B(\sqrt{2}u,\varepsilon
/\sqrt{2})}^{T}]\leq\mathcal{W}_{x_{0}}[|L_{B(u,\varepsilon)}^{T}|^{p/2}]\leq
CT^{\frac{p}{2}-1}\mathcal{W}_{x_{0}/\sqrt{2}}[L_{B(u/\sqrt{2},\sqrt
{2}\varepsilon)}^{T}] \label{dis2}%
\end{equation}
and (\ref{eq:Lmoments-bound4-bis}) is a consequence of the identity
\[
\mathcal{W}[L_{B(u,\varepsilon)}^{T}]=\varpi\varepsilon^{3}T.
\]
The proof is complete. \hfill$\square$

The previous lemma solves the main problem posed at the beginning of the
section. We prove also a related result in finite volume that we need for the
complementary results of section \ref{sectlocalizz}. We limit ourselves to the
upper bounds, for shortness.

\begin{lemma}
\label{lemmaultimo}Given $p>0$, $\alpha>0$, $u\in\mathbb{R}^{3}$, there is a
constant $C>0$ and a function $\theta\left(  R\right)  $ with $\lim
_{R\rightarrow\infty}\theta\left(  R\right)  =0$, such that, for $R>\left|
u\right|  +\ell$ and $T\geq\alpha\ell^{2}$%
\[
\int_{\left|  x_{0}\right|  \geq R}\mathcal{W}_{x_{0}}[|L_{B(u,\ell)}%
^{T}|^{p/2}]dx_{0}\leq C\ell^{p}\ell T\exp\left(  -\frac{R-\left(  \left|
u\right|  +\ell\right)  }{\sqrt{T}}\right)
\]
while for $T\leq\alpha\varepsilon^{2}$%
\[
\int_{\left|  x_{0}\right|  \geq R}\mathcal{W}_{x_{0}}[|L_{B(u,\varepsilon
)}^{T}|^{p/2}]dx_{0}\leq CT^{p/2}\varepsilon^{3}\exp\left(  -\frac{R-\left(
\left|  u\right|  +\varepsilon\right)  }{\sqrt{T}}\right)  .
\]
\end{lemma}

\textsc{Proof.}\hspace{0.3cm}Arguing as in the previous lemma, it is
sufficient to prove that
\[
\int_{\left|  x_{0}\right|  \geq R}\mathcal{W}_{x_{0}}[L_{B(u,\ell)}%
^{T}]dx_{0}\leq C\ell^{3}T\exp\left(  -\frac{R-\left(  \left|  u\right|
+\ell\right)  }{\sqrt{T}}\right)  .
\]
We have
\begin{align*}
\int_{\left|  x_{0}\right|  \geq R}\mathcal{W}_{x_{0}}[L_{B(u,\ell)}%
^{T}]dx_{0}  &  =\int_{B(0,R)^{c}}dx_{0}\int_{0}^{T}dt\int_{B(u,\ell
)}dz\;p_{t}\left(  z-x_{0}\right) \\
&  =\varpi\ell^{3}T\int_{0}^{T}\frac{dt}{T}\int_{B(u,\ell)}\frac{dz}%
{|B(u,\ell)|}\int_{B(0,R)^{c}}dx_{0}\;p_{t}\left(  z-x_{0}\right) \\
&  \leq\varpi\ell^{3}Tg\left(  R;T,u,\ell\right)
\end{align*}
with (recall that $\ell\leq1$)
\[
g\left(  R;T,u,\ell\right)  =\sup_{t\in(0,T],\left|  z\right|  \leq\left|
u\right|  +\ell}\int_{B(0,R)^{c}}dx_{0}\;p_{t}\left(  z-x_{0}\right)  .
\]
Now (denoting by $\left(  W_{t}\right)  $ a 3D Brownian motion)
\begin{align*}
g\left(  R;T,u,\ell\right)   &  =\sup_{t\in(0,T],\left|  z\right|  \leq\left|
u\right|  +\ell}P\left(  \left|  z-W_{t}\right|  \geq R\right) \\
&  \leq\sup_{t\in(0,T]}P\left(  \left|  W_{t}\right|  \geq R-\left(  \left|
u\right|  +\ell\right)  \right) \\
&  =\sup_{t\in(0,T]}P\left(  \left|  W_{1}\right|  \geq\frac{R-\left(  \left|
u\right|  +\ell\right)  }{\sqrt{t}}\right) \\
&  \leq P\left(  \left|  W_{1}\right|  \geq\frac{R-\left(  \left|  u\right|
+\ell\right)  }{\sqrt{T}}\right)
\end{align*}
This completes the proof. \hfill$\square$

%
%
%

\end{document}